\documentclass[a4paper,10pt]{article}
\usepackage[utf8x]{inputenc}

\usepackage{amsmath}
\usepackage{amssymb}
\usepackage{amsfonts}
\usepackage{amsthm,amscd}
\newtheorem{theorem}{Theorem}[section]
\newtheorem{lemma}[theorem]{Lemma}
\newtheorem{corollary}[theorem]{Corollary}
\newtheorem{proposition}[theorem]{Proposition}
\theoremstyle{definition}
\newtheorem{definition}[theorem]{Definition}

\theoremstyle{remark}
\newtheorem{remark}[theorem]{Remark}

\numberwithin{equation}{section}

\def\X{\mathbb X}

\title{\bfseries An enlargement of some symplectic objects\\ {\it Dedicated to Professor A. Banyaga
}}
\author{ \scshape St\'ephane Tchuiaga$^{}$\thanks{tchuiagas@gmail.com} }

\begin{document}

\maketitle

\begin{abstract}
The study of algebraic properties of groups of transformations
of a manifold gives rise to an interplay between different areas of mathematics such as topology, geometry, and  dynamical systems.
 Especially, in this paper, we point out some  
interplays between  topology, geometry, and dynamical systems which are underlying to the group of symplectic homeomorphisms. 
The latter situation can occur when one thinks of the following question : Is there a $C^0-$ flux geometry 
which is underlying to the group of strong symplectic homeomorphisms so that Fathi's Poincar\'e duality theorem continues to hold? 
In this paper, we discuss on some possible answers of the above preoccupation, and we elaborate various topological analogues of some well-known
 results found in the field of symplectic dynamics.\\

{\bf AMS Subject Classification:} 53D05, 53D35, 57R52, 53C21.\\
{\bf Key Words and Phrases:} 
Symplectic displacement energy, Hofer-like norm, Hodge decomposition theorem, Symplectic isotopies, Rigidity,
 Homeomorphisms, Diffeomorphisms, Symplectic homeomorphisms, Homotopy, 
Topological symplectic isotopies, de Rham Cohomology, Poincar\'e duality theorem, Fathi's mass flow.

\end{abstract}
\begin{center}
 \section{Introduction}
\end{center}
 
According to Oh-M\"{u}ller \cite{Oh-M07}, the automorphism group of the $C^0-$symplectic topology is the closure of the 
group $Symp(M,\omega)$ of all symplectomorphisms of a symplectic manifold  $ (M,\omega) $
in the group  $Homeo(M)$ of all homeomorphisms of $M$ equipped with the $C^0-$topology. That group, denoted $Sympeo(M,\omega)$ has 
been called group of all symplectic homeomorphisms :  
$$Sympeo(M,\omega) = \overline{Symp(M,\omega)}\subset Homeo(M).$$ 
This definition has been motivated by the following celebrated 
rigidity theorem dues to Eliashberg \cite{Elias87} and Gromov \cite{MGrom85}.\\
\begin{theorem}(\cite{Elias87, MGrom85})\label{E-G}
 The group $Symp(M,\omega)$ of all symplectomorphisms of a symplectic manifold  $ (M,\omega) $ is $C^0$ closed 
inside the group of diffeomorphisms over $M$. 
\end{theorem}
Oh-M\"{u}ller showed that any symplectic homeomorphism preserves the Liouville measure \cite{Oh-M07}.
 Furthermore, a result dues to Oh-M\"{u}ller states 
that the group of Hamiltonian homeomorphisms is contained in the kernel of Fathi's mass flow. 
This group is a proper subgroup in the identity component of $Sympeo(M,\omega)$. Buhovsky \cite{Buh}  
observed that Eliashberg-Gromov $C^0-$rigidity follows
from Oh-M\"uller theorem on the uniqueness of topological Hamiltonians of topological Hamiltonian systems \cite{Oh-M07}. All 
the above facts demonstrate once more an interest of the study of symplectic homeomorphisms.\\ 

Recently, 
 motivated by the result in Theorem \ref{E-G}, Banyaga  \cite{ Ban08a, Ban10c} defined two classes 
of symplectic homeomorphisms called the {\bf  strong symplectic homeomorphisms in the $L^{\infty}-$context}  
and  the {\bf  strong symplectic homeomorphisms in the $L^{(1,\infty)}-$context}. 
The two contexts of { \bf strong symplectic homeomorphisms} arise from two different topologies, but  
 it is proved in Banyaga-Tchuiaga \cite{BanTchu} that the nature of any { \bf strong symplectic homeomorphism}
 does not depend on the choice of the  $L^{\infty}$ symplectic topology 
 or  the $L^{(1,\infty)}$ symplectic topology. 
To understand the statement of the mains result of this paper, 
we will need the following definition from \cite{BanTchui02, Tchuia03}. 

\begin{definition}(\cite{BanTchui02, Tchuia03})\label{SSI}
 A continuous family $(\gamma_t)$ of homeomorphisms of $M$ with 
$\gamma_0 = id$ is called a strong symplectic isotopy (or ssympeotopy ) if there exists a sequence $\Phi_i$ of symplectic isotopies 
which converges uniformly to 
$(\gamma_t)$ such that the sequence of symplectic vector fields generated by $\Phi_i$ is Cauchy in the $L^\infty-$ 
Hofer-like norm defined in \cite{Ban08a}.
\end{definition}
If the manifold is simple connected, then any ssympeotopy is called a continuous Hamiltonian flow ( or Hameotopy) in the sense of 
Oh-M\"{u}ller.\\

A result from  
 Banyaga-Tchuiaga \cite{BanTchui02}
shows that in the above definition, if we set $\Phi_i = (\phi_i^t) $ and consider $Z_t^i(x) = \dfrac{d\phi_i^t}{dt}\circ(\phi_i^t)^{-1}(x),$ 
for each $i$ and for all $x\in M$, then 
the limit $\mathcal{M}(\gamma)$ of the sequence of symplectic vector fields $(Z_t^i)$ w.r.t 
the $L^\infty-$ 
Hofer-like norm is independent of the choice of 
the sequence of symplectic isotopies $\Phi_i$ in Definition \ref{SSI}. The authors called the latter limit the ''generator'' of the
ssympeotopy $(\gamma_t)$. It follows from Banyaga-Tchuiaga \cite{BanTchui02} that any generator of strong symplectic isotopy 
is of the form $(U,\mathcal{H})$ where the maps $(t,x)\mapsto U_t(x)$ and $t\mapsto\mathcal{H}_t$ are continuous, and for each $t$, $\mathcal{H}_t$ 
is a smooth harmonic $1-$form for some Riemannian metric $g$ on $M$. Here is the main result from \cite{BanTchui02}.\\
\begin{theorem}\label{URBT}(\cite{BanTchui02})
 Any strong symplectic isotopy  determines a unique generator. 
\end{theorem}
This theorem generalizes the uniqueness theorem of generating functions for continuous Hamiltonian flows from Viterbo \cite{Viter06, Viter92}, and
Buhovsky-Seyfaddini \cite{Buh-Fay10}. Using the above uniqueness result of generators of ssympeotopies, Banyaga pointed 
out that any smooth hameomorphisms is a Hamiltonian diffeomorphisms, a theorem on Hamiltonian rigidity. 
\\ However, in the presence of a positive symplectic displacement energy from 
Banyaga-Hurtubise-Spaeth \cite{BDS}, we point out the following converse of Theorem \ref{URBT}.

\begin{theorem}\label{UGG}
 Any generator of strong symplectic isotopy corresponds to a unique strong symplectic isotopy.
\end{theorem}
  The above result    
generalizes a result found by Oh (Theorem $3.1$, \cite{Oh}) in the study of topological Hamiltonian dynamics. 
Following the proof of  Theorem $1.2$ found in \cite{TD}, it seems that another proof of Theorem \ref{UGG} can be pointed out  
without appealing to the positivity of the symplectic displacement energy. But, this other way to prove Theorem \ref{UGG} is a little longer than the proof  
involving the symplectic displacement energy. To show an interest of the uniqueness results of Theorem \ref{URBT} and 
 Theorem \ref{UGG} in the study of topological symplectic dynamics, let's go back to the smooth symplectic geometry.\\
 
Using Hodge's decomposition theorem of differential forms, one can establish that there is a one-to-one 
correspondence between the space of all symplectic isotopies and the space of 
 all the pairs $(V,\mathcal{K})$ where $V$ is a smooth family of normalized function and $\mathcal{K}$ 
is a smooth family of smooth harmonic $1-$form for some Riemannian metric $g$ on $M$ 
(any differentiable manifold $M$ can be equipped with a Riemannian metric). 
The latter pairs was called the smooth generators of symplectic paths (see \cite{BanTchu}). A result from \cite{BanTchu} states that there 
is a group isomorphism between the space of symplectic paths and that of all their generators. 
This suggests that the flux 
of any symplectic isotopy $\Phi$ generated by $(V,\mathcal{K})$ is exactly the de Rham cohomology class $ \int_0^1 [\mathcal{K}_t]dt$  
where $[.]$ stands for the de Rham cohomology class. Visibly, this way to evaluate the flux of any symplectic isotopy does 
not involves the derivative of the corresponding isotopy.\\ This seems to suggest that in the presence 
of the uniqueness results of Theorem \ref{URBT} and 
 Theorem \ref{UGG}, we can assign to any ssympeotopy $(\gamma_t)$ generated by $(U,\mathcal{H})$ a well defined 
de Rham cohomology class $ \int_0^1 [\mathcal{H}_t]dt$ since the map $t\mapsto\mathcal{H}_t$ is continuous, 
and for each $t$, $\mathcal{H}_t$ 
is a smooth harmonic $1-$form. This will be clearly exposed in greater detail later on.\\ More generally,
 it is known that on any closed symplectic manifold, 
various constructions involving symplectic paths (e.g, Hofer metrics, Hofer-like metrics, symplectic flux etc..) are based 
on the one-to-one correspondence that exists between the space of such paths and that of smooth 
families of smooth symplectic vector fields. Actually, the uniqueness results  of Theorem \ref{URBT} and 
 Theorem \ref{UGG} tell us to investigate the following questions: 
\begin{itemize}
\item Is there a $C^0-$flux geometry 
which is underlying to the group of { \bf strong symplectic homeomorphisms}? 
\item If the answer of the above question is affirmative, then 
can Fathi's Poincar\'e duality theorem continue to hold?
\item Is Banyaga's Hofer-like geometry admits a topological analogue 
in the world of { \bf strong symplectic homeomorphisms}? 
\item Since the mass flow of any { \bf strong symplectic isotopy} exists, can we  
give an explicit formula of the latter mass flow?

\end{itemize}

The goal of this paper is to investigate the above questions. \\

We organize the present paper as follows. 
In Section 2, we recall some fundamental facts concerning symplectic mappings and isotopies. Section 2.1, deals with a
description of symplectic isotopies that was introduced in \cite{BanTchu}.
 In Sections 2.2, 2.3, and 2.4 we 
recall the definitions of Banyaga's Hofer-like topologies and 
symplectic displacement energy.  Section 3 deals 
with Hodge-de Rham’s theory, flux for symplectic maps, and Fathi's mass flow.
Here, we recall the constructions of flux for symplectic maps and Fathi's mass flow; 
including de Rham’s theorem, the Poincar\'e duality theorem, and Fathi's Poincar\'e duality theorem. 
In Section 4, we will recall the definitions of topological symplectic isotopies.
 Moreover, we show an impact of Fathi's Poincar\'e duality theorem 
 in the study of some homotopic property (relatively to fixed extremities) 
of certain class of strong symplectic isotopies. As a consequence of the latter study, we derive that 
 any ssympeotopies whose Fathi's mass flow 
is trivial admits its extremities in the group of Hamiltonian homeomorphisms. Once more, using 
a result that was proved by Weinstein \cite{AWein}, we prove a result which implies in turn that 
the flux group $\Gamma$ stays 
invariant under a perturbation of the fundamental group $\pi_1(Symp_0(M,\omega))$ withe respect to a suitable $C^0-$symplectic topology. We 
shall 
also prove that any strong symplectic isotopy which is a $1-$parameter group  
decomposes as composition of a smooth harmonic flow and a continuous Hamiltonian flow in the sense of Oh-M\"{u}ller. 
Section 5, we will introduces an enlargement of the symplectic flux homomorphism. Here, we give an explicit 
formula for computing the Fathi mass flow of any strong symplectic isotopy, and construct a surjective 
group homomorphism from the group of strong symplectic isotopies onto the quotient of the first de Rham cohomology group by the flux group $\Gamma$. 
We prove that its kernel coincides with the group of Hamiltonian homeomorphisms. 
This map is the $C^0-$analogue of 
the well known flux homomorphism.  As a consequence of the latter fact, we 
see that the group of Hamiltonian homeomorphisms is path connected and locally connected, while the group 
of strong symplectic homeomorphism is locally path connected. In Section 6, we 
construct a $C^0-$Hofer-like norm for
 strong symplectic homeomorphisms, and prove that the restriction of the latter norm to the group 
of Hamiltonian homeomorphisms is equivalent to the norm constructed by Oh on the space of 
Hamiltonian homeomorphisms. Section 7 deals with an enlargement of the symplectic displacement energy. 
Here, we use the result of Section 6 to define a positive symplectic displacement energy for strong symplectic 
homeomorphisms. Finally, Section 8 gives some examples and introduce some conjectures.

\begin{center}
 \section{Preliminaries}
\end{center}
Let $M$ be a smooth closed manifold of dimension $2n$. A differential $2-$form $\omega$ on $M$ is called a symplectic form if
 $\omega$ is closed and nondegenerate. In particular, any symplectic manifold is oriented.  
From now on, we shall always assume that $M$ admits a symplectic form $\omega$. A diffeomorphism $\phi:M\rightarrow M$ 
is called symplectic if it preserves the symplectic form $\omega$, i.e. $\phi^\ast(\omega) = \omega.$ We denote 
by $Symp(M,\omega)$ the symplectomorphisms' group.\\

\subsubsection{Symplectic vector fields}
The symplectic structure $\omega$ on $M,$ being nondegenerate, induces an isomorphism between vector fields $Z$ and $1-$forms on $M$ given by 
$Z\mapsto \omega(Z,.)=: \iota(Z)\omega$. 
A vector field $Z$ on $M$ is symplectic if $\iota(Z)\omega$ is closed. In particular, 
a symplectic vector field $Z$ on $M$ is said to be  a Hamiltonian vector 
field if $\iota(Z)\omega$ is exact. It follows from the definition of symplectic vector fields that, 
 if the first de Rham cohomology group of the manifold $M$ is trivial (i.e. $H^1(M,\mathbb{R}) = 0$), 
then all the symplectic vector fields induced  
by a symplectic form $\omega$ on $M$ are Hamiltonian. If we equip $M$ with a Riemannian metric $g$ 
(any differentiable manifold $M$ can be equipped with a Riemannian metric), then any harmonic $1-$form $\alpha$ on $M$ 
determines a symplectic
 vector field $Z$ such that $\iota(Z)\omega = \alpha$ (so-called harmonic vector field, see \cite{Ban08a}).\\ 
\subsubsection{Symplectic isotopies}
An isotopy $\{\phi_t\}$  of the symplectic manifold $(M,\omega)$ is said to be symplectic if for each $t$, the vector field 
$Z_t = \dfrac{d\phi_t}{dt}\circ\phi_t^{-1}$ is symplectic. In particular 
 a symplectic isotopy $ \{\psi_t\}$ is called Hamiltonian if for each $t$,  
$Z_t = \dfrac{d\psi_t}{dt}\circ\psi_t^{-1}$ is Hamiltonian, i.e. there exists a smooth 
function $F: [0,1]\times M \rightarrow \mathbb{R}$ called Hamiltonian such that
 $\iota(Z_t)\omega = dF_t$. As we can see, any Hamiltonian isotopy determines a 
 Hamiltonian $F: [0,1]\times M \rightarrow \mathbb{R}$ up to an additive constant . 
Throughout the paper we assume that all Hamiltonians are normalized in the following way: given a Hamiltonian $F: [0,1]\times M \rightarrow \mathbb{R}$
we require that $\int_M F_t \omega^n = 0$. We denote by $\mathcal{N}([0,1]\times M\ ,\mathbb{R})$
the space of all smooth normalized Hamiltonians and by $Ham(M,\omega)$ 
the set of all time-one maps of Hamiltonian isotopies. If we equip $M$ with a Riemannian metric 
$g$, then a symplectic isotopy $ \{\theta_t\}$ is said to be harmonic if for each $t$,
 $Z_t = \dfrac{d\theta_t}{dt}\circ\theta_t^{-1}$ is harmonic. We denote 
by $Iso(M, \omega)$ the group of all symplectic isotopies of $(M,\omega)$ and by $Symp_0(M,\omega)$ 
the set of all time-one maps of symplectic isotopies.

\subsubsection{Harmonics $1-$forms}
 From now on, we assume that $M$ is equipped with a Riemannian metric $g$, 
and denote by $\mathcal{H}^1(M,g)$ the space of harmonic $1-$forms on $M$ 
with respect to the Riemannian metric $g$. In view of the Hodge theory,  $\mathcal{H}^1(M,g)$ 
is a finite dimensional vector space over $\mathbb{R}$ which is isomorphic to $H^1(M,\mathbb{R})$ (see \cite{FW71}). 
The dimension of $\mathcal{H}^1(M,g)$ is the first Betti number of the manifold $M$, denoted by 
$b_1$. Taking ${(h^i)_{1\leq i\leq b_1}}$ as a basis of the vector space $\mathcal{H}^1(M,g)$, we 
 equip $\mathcal{H}^1(M,g)$ with the Euclidean norm $|.|$ defined as follows : 
for all $H$ in $\mathcal{H}^1(M,g)$ with $H = \Sigma_{i = 1}^{b_1}\lambda^ih^i$ 
we have $|H| :=\Sigma_{i = 1}^{b_1}|\lambda^i|.$ It is convenient to compare the above Euclidean norm with the well-known 
uniform sup norm of differential $1-$forms. For this purpose, let's recall the definition of the uniform sup norm of a differential $1-$form
 $\alpha$ on $M$. For all $x\in M$, we know that $\alpha$ induces a linear map $\alpha_x : T_xM \rightarrow \mathbb{R}$ 
whose norm is given by $\|\alpha_x\| = \sup\{ |\alpha_x(X)| : X\in T_xM,  \|X\|_g = 1\},$
where $\|.\|_g$ is the norm induced 
on each tangent space $ T_xM$ by the Riemannian metric $g$. Therefore, the uniform sup norm of $\alpha$, say $|.|_{0}$ is defined by 
 $|\alpha|_{0} = \sup_{x\in M}\|\alpha_x\|.$
 In particular, when $\alpha$ is a harmonic $1-$form (i.e $\alpha = \Sigma_{i = 1}^{b_1}\lambda^ih^i$), we obtain the following estimates :  
$|\alpha|_{0} \leq \Sigma_{i = 1}^{b_1}|\lambda^i| | h^i|_0\leq E|\alpha|$
where $ E := \max_{1\leq i\leq b_1}|h^i|_{0}$. 
If the basis ${(h^i)_{1\leq i\leq b_1}}$ is such that $E \textgreater 1$, then one can always normalize such a basis so that  
 $E$ equals $1$. Otherwise, the identity $|\alpha|_{0}\leq E|\alpha|$ reduces to $|\alpha|_{0}\leq|\alpha|$.
We denote by $ \mathcal P\mathcal{H}^1(M,g)$, the space of smooth mappings $\mathcal{H}: [0,1]\rightarrow\mathcal{H}^1(M,g)$.

 \subsection{A description of symplectic isotopies \cite{BanTchu}}\label{SC3}
In this subsection, from the group of symplectic isotopies, we shall deduce another group which
will be convenient later on (see \cite{BanTchu}). Consider $\{\phi_t\}$ to be a symplectic isotopy, for each $t$, the vector field 
$Z_t = \dfrac{d\phi_t}{dt}\circ(\phi_t)^{-1}$ satisfies $d\iota(Z_t)\omega =0$. So, it follows from   
 Hodge's theory that
$\iota(Z_t)\omega$ decomposes as the sum of an exact $1-$form $dU_t^\Phi$ and a harmonic $1-$form $\mathcal{H}_t^\Phi$ (see \cite{FW71}). 
Denote by $U$ the Hamiltonian $U^\Phi = (U^\Phi_t)$ normalized,  
and by $ \mathcal{H}$ the smooth family of harmonic $1-$forms $ \mathcal{H}^\Phi = (\mathcal{H}^\Phi_t)$. 
In \cite{BanTchu}, the authors denoted by $\mathfrak{T}(M, \omega, g)$ the Cartesian product 
$\mathcal{N}([0,1]\times M ,\mathbb{R})\times \mathcal P\mathcal{H}^1(M,g)$, and equipped it with a
 group structure which makes the bijection
\begin{equation}
  Iso(M, \omega)\rightarrow\mathfrak{T}(M, \omega, g), \Phi\mapsto (U, \mathcal{H})
\end{equation}
 a group isomorphism. Denoting the map just constructed by $\mathfrak{A}$, the authors denoted  
any symplectic isotopy $\{\phi_t\}$ as $\phi_{(U,\mathcal{H})}$ 
to mean that the mapping $\mathfrak{A}$ maps $\{\phi_t\}$ onto $(U,\mathcal{H})$, and $(U,\mathcal{H})$ is 
called the 
``generator'' of the symplectic path $\phi_{(U,\mathcal{H})}$. In particular, any symplectic isotopy of the form 
$\phi_{(0,\mathcal{H})}$ is considered to be a harmonic isotopy, while any 
symplectic isotopy of the form $\phi_{(U,0)}$ is considered to be a Hamiltonian 
isotopy. The product in $\mathfrak{T}(M, \omega, g)$ is given by, 
\begin{equation}
  (U,\mathcal{H})\Join(V ,\mathcal{K}) = ( U + V\circ\phi_{(U,\mathcal{H})}^{-1} 
+ \widetilde \Delta(\mathcal{K},\phi_{(U,\mathcal{H})}^{-1}), \mathcal{H} + \mathcal{K} )
\end{equation}
The inverse of $(U,\mathcal{H})$, denoted $\overline{(U,\mathcal{H})}$ is given by 
\begin{equation}
 \overline{(U,\mathcal{H})} = (- U\circ\phi_{(U,\mathcal{H})} - \widetilde\Delta(\mathcal{H},\phi_{(U,\mathcal{H})}),-\mathcal{H})
\end{equation}
where for each $t$, $\phi_{(U,\mathcal{H})}^{-t} :=(\phi_{(U,\mathcal{H})}^{t})^{-1}$, and 
$\widetilde \Delta_t(\mathcal{K},\phi_{(U,\mathcal{H})}^{-1})$ is the function\\ 
$$\Delta_t(\mathcal{K},\phi_{(U,\mathcal{H})}^{-1}) 
:= \int_0^t\mathcal{K}_t(\dot{\phi}_{(U,\mathcal{H})}^{-s})\circ\phi_{(U,\mathcal{H})}^{-s}ds,$$ normalized.

\subsection{Banyaga's Topologies \cite{BanTchu}, \cite{Tchuia03}}\label{SC4}
In this subsection, we reformulate the topology introduced by Banyaga \cite{Ban08a} on 
the space of symplectic vector fields.  Let $X$ be a symplectic vector field. The $1-$form 
$i_X\omega$ can be decomposed in a unique way as the sum of a harmonic $1-$form $ \mathcal{H}_X$ with an exact $1-$form $dU_X$.  
The function $U_X$ is given by $U_X = \delta G(i_X\omega),$ where $\delta$ is the codifferential 
operator and $G$ is the Green operator \cite{FW71}.   
In regard of the above decomposition of symplectic vector fields, we will denote any symplectic 
vector field $X$ by $X^{(U,\mathcal{H})}$ to mean that Hodge's decomposition of $\iota(X)\omega$
 gives $dU + \mathcal{H}$ with $U$ normalized. According to Banyaga, the above decomposition of 
symplectic vectors  gives rise to an intrinsic norm on the space 
of symplectic vector fields, defined by $$\|X^{(U,\mathcal{H})}\| = |\mathcal{H}| + osc(U),$$ where 
$osc(U) = \max_{x\in M}U(x) - \min_{x\in M}U(x),$ and $|.|$ represents the Euclidean norm on the 
space of harmonic $1-$forms that we introduced in the beginning.
The norm $\|.\|$ gives rise to 
a norm defined on the space of $1-$parameter symplectic vectors fields as follows. Let $(Y_t)_t$ be a smooth family of 
symplectic vector fields, we have : 
$$\|(Y_t)\|^\infty = \max_t\|Y_t\|.$$ 
The above norm is called the $L^\infty-$ Banyaga's norm of the family of symplectic vector fields $(Y_t)_t.$ 
For instance, the above norm induces a distance on the space $\mathfrak{T}(M, \omega, g)$ as follows :  For all 
$(U,\mathcal{H}), (V ,\mathcal{K})\in\mathfrak{T}(M, \omega, g)$, 
\begin{equation}\label{metric2}
  D^2 ((U,\mathcal{H}), (V ,\mathcal{K})) = (\|(X^{(U_t,\mathcal{H}_t)} - X^{(V_t,\mathcal{K}_t)})_t\|^\infty + 
\|(X^{\overline{(U_t,\mathcal{H}_t)}} - X^{\overline{(V_t ,\mathcal{K}_t)}})_t\|^\infty)/2,
\end{equation} 
 Therefore, the $L^\infty-$topology on the space  $\mathfrak{T}(M, \omega, g)$  is the one induced by the metric $D^2$.

\subsection{Banyaga's Hofer-like norms \cite{Ban08a}}\label{SC5} 
According to \cite{Ban08a}, 
the so-called $L^{(1,\infty)}$ version and $L^{\infty}$ version of Banyaga's Hofer-like lengths of any  
$\Phi = \phi_{(U,\mathcal{H})}\in Iso(M,\omega)$ are defined respectively by, 
\begin{equation}\label{blg1}
 l^{(1,\infty)}(\Phi) = \int_0^1 osc(U_t) + |\mathcal{H}_t|dt, 
\end{equation}
\begin{equation}\label{blg2}
 l^\infty(\Phi) 
= \max_{t\in[0,1]}( osc(U_t) + |\mathcal{H}_t|).
\end{equation}
Clearly $l^{(1,\infty)}(\Phi) \neq l^{(1,\infty)}(\Phi^{-1})$ unless $\Phi$ is Hamiltonian. Indeed,  $\Phi = \phi_{(U ,\mathcal{H})}$ 
implies that  $\Phi^{-1}
= \phi_{\overline{(U,\mathcal{H})}}$ where $\overline{(U,\mathcal{H})} = 
(-U\circ\Phi - \widetilde\Delta( \mathcal{H},\Phi),-\mathcal{H})$. 
Hence, we see that the mean oscillation of the function $U$ can be different from that of the function
$-U\circ\Phi - \widetilde\Delta( \mathcal{H},\Phi)$. But, if $\Phi$ is Hamiltonian, i.e. 
$\Phi = \phi_{(U , 0)}$, then the mean oscillation of $ U $ is equal to that
of $-U\circ\Phi$, i.e. $l^{(1,\infty)}(\Phi) = l^{(1,\infty)}(\Phi^{-1})$.  
Similarly, we have $l^{\infty}(\Phi) \neq l^{\infty}(\Phi^{-1})$ unless $\Phi$ is Hamiltonian.  
Now, let $\phi\in Symp_0(M,\omega)$, using the above Banyaga's lengths, 
Banyaga \cite{Ban08a} defined respectively
the  $L^{(1,\infty)}$ energy and $L^{\infty}$ energy of $\phi$ by,
\begin{equation}\label{bener1}
 e_0(\phi) = \inf(l^{(1,\infty)}(\Phi)),
\end{equation}
\begin{equation}\label{bener2}
 e^{\infty}_0(\phi) = \inf(l^\infty(\Phi)),
\end{equation}
where the infimum are taken over all symplectic isotopies $\Phi$ with time-one map equal to $\phi$. 
Therefore,  
the $L^{(1,\infty)}$ Banyaga's Hofer-like norm and the $L^{\infty}$ Banyaga's Hofer-like norm of $\phi$ are respectively defined by,
\begin{equation}\label{bny1}
 \|\phi\|_{HL}^{(1,\infty)} = (e_0(\phi) + e_0(\phi^{-1}))/2,
\end{equation}
\begin{equation}\label{bny2}
 \|\phi\|_{HL}^\infty = (e^\infty _0(\phi) + e^\infty _0(\phi^{-1}))/2.
\end{equation}
Each of the norms $ \|.\|_{HL}^\infty $ and $\|.\|_{HL}^{(1,\infty)}$ 
generalizes the Hofer norms for Hamiltonian diffeomorphisms in the following sense : 
in the special case of a closed symplectic manifold $(M,\omega)$ for which $Ham(M,\omega) = Symp_0(M,\omega)$ (or  $H^1(M,\mathbb{R}) = 0$),
 the norm $ \|.\|_{HL}^\infty $ reduces to a norm $ \|.\|_{H}^\infty $ called the $L^\infty$ Hofer norm, while 
the norm $\|.\|_{HL}^{(1,\infty)}$ reduces to a norm $ \|.\|_{H}^{(1,\infty)} $ called the $L^{(1,\infty)}$ Hofer
 norm. But, a result that was proved by Polterovich \cite{Polt93} shows that the above
 Hofer's norms are equal in general, i.e. $ \|.\|_{H}^{(1,\infty)}=  \|.\|_{H}^{\infty} $. In other words,  
the norms $ \|.\|_{HL}^\infty $ and $\|.\|_{HL}^{(1,\infty)}$ are equal when $Ham(M,\omega) = Symp_0(M,\omega)$ (or  $H^1(M,\mathbb{R}) = 0$).
It was  proved in \cite{TD} that the following equality holds true  
$ \|.\|_{HL}^\infty = \|.\|_{HL} .$ This points out the uniqueness of Banyaga's Hofer-like geometry, and then generalizes a result 
that Polterovich \cite{Polt93} proved (Lemma 5.1.C, \cite{Polt93}).
\subsection{Displacement energy (Banyaga-Hurtubise-Spaeth)}
\begin{definition}(\cite{BDS})
 The symplectic displacement energy $ e_S(A)$ of a non empty set $A\subset M$ is :
$$e_S(A) = \inf\{\|\phi\|_{HL} |  \phi\in Symp(M,\omega)_0, \phi(A)\cap A = \emptyset\}.$$
\end{definition}
\begin{theorem}\label{BDS}(\cite{BDS})
 For any non empty open set $A\subset M$, $ e_S(A)$ is a strict positive number.
\end{theorem}

A straightforward application of the equality $\|,\|_{HL} = \|,\|_{HL}^\infty$ implies
 that  the above definition of positive symplectic displacement energy does not depend on the choice of Banyaga's Hofer-like norm. 
Furthermore, in the Section 5 of the present paper, we will see how one extends the above displacement symplectic energy in the world of 
strong symplectic homeomorphisms. This will be  
supported by the uniqueness result from \cite{BanTchu} and the uniqueness of Banyaga's Hofer-like geometry \cite{TD}. 

%
%

\section{Hodge-de Rham’s theory}
\subsection{De Rham’s theorem}
We denote by $\Omega^k(M)$ the space of smooth differential $k−$forms on $M$. 
In local coordinates $(x_1 , x_2 , ...,x_{2n} )$, an element $\alpha \in \Omega^k(M)$ has the following 
expression $$\alpha = \Sigma_{ I = \{i_1\textless i_2\textless, ...,\textless i_k\}}a_Idx_1\wedge dx_2.....\wedge dx_k = \Sigma_I a_I dx_I,$$
where $a_I$ is a smooth function of $(x_1 , x_2 , ..., x_{k} )$. The exterior differentiation is a
differential operator $$ d : \Omega^k(M)\rightarrow C^\infty(\wedge^{k +1} T^\ast M ).$$
Locally, $$d(\Sigma_I a_I dx_I) = \Sigma_I da_I \wedge  dx_I.$$ This operator satisfies $d \circ d = 0,$ hence the range of $d$ is included in the kernel
of $d$. 
\begin{definition}(\cite{FW71}) 
 The $k^{th}$ de Rham’s cohomology group of $M$ is defined by
$$H^k_{dH}(M) = \dfrac{\{\alpha \in \Omega^k(M), d \alpha = 0\}}{ d (\Omega^{k-1}(M))}.$$
\end{definition}
These spaces are clearly diffeomorphism invariants of M. Moreover a deep
theorem of G. de Rham says that these spaces are isomorphic to the real cohomology
group of M. Thus, the are homotopic invariants.

\subsection{The flux (Banyaga)}
Denote by $Diff^\infty(M)$ the group of diffeomorphisms on $M$ endowed with the $C^\infty$ compact open topology. Let $Diff_0(M)$ be the 
identity component of $Diff^\infty(M)$ for 
the $C^\infty$ compact open topology. Let $\Omega$ be any closed $p-$form on $M$.
Denote by $Diff^\Omega(M)\subset Diff^\infty(M)$ the space of all diffeomorphisms that preserve the $p-$form $\Omega$,
 and by $Diff_0^\Omega(M)$ we denote the connected component by smooth arcs of the identity in $Diff^\Omega(M)$. 
Let $\Phi = (\phi^t)$ be a smooth path in $Diff_0^\Omega(M)$ with $\phi_0 = Id$, and set  
\begin{equation}
 \dot{\phi}_t(x) = \dfrac{d\phi_t}{dt}(\phi_t^{-1}(x))
\end{equation}
for all $(t, x)\in [0,1]\times M.$ It was proved in \cite{Ban78} that 
\begin{equation}
 \sum(\Phi) = \int_0^1(\phi_t^\ast(i_{\dot{\phi}_t}\Omega))dt,
\end{equation}
is a closed $p-1$ form and its cohomology class denoted by $[\sum(\Phi)]\in  H^{(p-1)}(M,\mathbb{R})$ depends only on the homotopy
class $\{\Phi\}$ of the isotopy $\Phi$ relatively with fixed ends in  $Diff^\Omega_0(M)$, and the map 
$\{\Phi\}\mapsto [\sum(\Phi)]$ is a surjective group homomorphism
$$Flux_\Omega : \widetilde{Diff_0^\Omega(M)}\rightarrow  H^{(p-1)}(M,\mathbb{R}).$$ 

In case $\Omega$ is a symplectic form $\omega$, we get a homomorphism 
$$Flux_\omega : \widetilde{Symp_0(M,\omega)}\rightarrow {H}^1(M, \mathbb{R}),$$ 
where $\widetilde{Symp_0(M,\omega)}$ is the universal covering of the space $Symp_0(M,\omega)$.\\
 Denote by $\Gamma$ the image by $Flux_\omega$ of 
$\pi_1(Symp_0(M,\omega))$. The homomorphism 
$ Flux_\omega$ induces a surjective homomorphism $flux_\omega$ 
 from $Symp_0(M,\omega)$ onto $ {H}^1(M, \mathbb{R})/\Gamma.$ 
From the above construction, Banyaga \cite{Ban78, Ban97} proved that the group of all Hamiltonian diffeomorphisms of any compact 
symplectic manifold $(M,\omega)$ is a simple group which coincides with the kernel of $flux_\omega$, a very deep result. We will 
need the following result from \cite{Ban78, Ban97}. 
\begin{theorem}\label{Fgeo} Let $\Phi$ be a symplectic isotopy. 
\begin{enumerate}
\item
If $|Flux_\omega(\Phi)|$ is arbitrarily small, then 
 any symplectic isotopy $\Psi$ with the same extremities than $\Phi$ has its flux in $ \Gamma$, i.e. 
$Flux_\omega(\Psi)\in \Gamma$. 

 \item If  $Flux_\omega(\Phi) = 0$, then $\Phi$ is homotopic with fixed endpoints to a Hamiltonian isotopy. 
\end{enumerate}
\end{theorem}
{\it Proof}. Let $\Phi = (\phi^t) $ and $\Psi = (\psi^t)$ be two symplectic isotopies with the same endpoints such that
 $|Flux_\omega(\phi^t)|$ is arbitrarily small. Since $ \Psi\circ\Phi^{-1}$ is a loop at the identity, 
then 
 $Flux_\omega(\Psi\circ\Phi^{-1})\in \Gamma$. But, 
$$d(Flux_\omega(\Psi), \Gamma) = \inf_{\gamma\in \Gamma} |\gamma - Flux_\omega(\Psi)|$$ 
$$\leq | Flux_\omega(\Psi\circ\Phi^{-1}) - Flux_\omega(\Psi)| $$ 
$$\leq |  Flux_\omega(\Phi)|,$$
where the right hand side is arbitrarily small. Thus, $Flux_\omega(\Psi)\in \Gamma$ 
since $\Gamma$ is discrete (see Ono \cite{Ono}). For (2), we follow the proof of 
a result from MacDuff-Salamon \cite{MacDuf-SAl} (Theorem 10.12, \cite{MacDuf-SAl}).
Assume that $\Phi = (\phi_t)$ is generated by $(U,\mathcal{H})$. By assumption we have 
$$Flux_\omega(\Phi) =  Flux_\omega(\phi_{(0,\mathcal{H})}^t) = [\int_0^1\mathcal{H}_tdt] = 0.$$ Thus, it follows from 
Hodge's theory that $\int_0^1\mathcal{H}_tdt = 0$ since $\int_0^1\mathcal{H}_tdt$ is harmonic and $M$ is compact. 
For all $t$, as in \cite{MacDuf-SAl}, let $\theta^t_s$ be the flow generated by the symplectic 
vector field $Y_t = -(\int_0^t\mathcal{H}_udu)^\sharp,$  i.e. $Y_t = \dfrac{d\theta_s^t}{ds}\circ(\theta_s^t)^{-1}$ with $s\in \mathbb{R},$ 
and $\theta^t_0 = id_M$. 
Since $Y_1 = 0 = Y_0$, we derive that $ \theta^0_s = id_M = \theta^1_s$ for all $s$. 
It is not too hard to see that the isotopy $\varphi_t = \theta^t_1\circ\phi_{(0,\mathcal{H})}^t $ 
is homotopic to $ \phi_{(0,\mathcal{H})}^t$ relatively with fixed endpoints. In addition, 
  for each $u\in[0,1]$, we have $Flux_\omega((\varphi_t)_{0\leq t\leq u}) =0$. This implies that 
$$\int_0^u[\iota(\dot\varphi_t)\omega ]dt = 0,$$ for all $u\in[0,1]$, i.e. $\iota(\dot\varphi_t)\omega$ is exact for all $t\in[0,1]$. 
Thus, the isotopy $(\varphi_t)_t$ is Hamiltonian with $\varphi_1 =\phi_{(0,\mathcal{H})}^1.$ Let $(\alpha_t)$ 
be the Hamiltonian part in Hodge's decomposition of $\Phi$.  
The mapping $$ H : [0, 1]\times [0, 1]\rightarrow Symp_0(M,\omega),$$ $$
 (s,t)\mapsto \theta^t_s\circ \phi_t,$$ induces a homotopy between $(\varphi_t\circ\alpha_t)_t$ and $\Phi$. $\Box$

\subsection{The mass flow (Fathi, \cite{Fathi80} )}
Let $\mu$ be a ''good measure'' on the manifold $M$. Let  $Homeo_0(M,\mu)$ 
 denotes the identity component 
in the group of measure preserving homeomorphisms $Homeo(M,\mu)$, and 
  $\widetilde {Homeo_0(M,\mu)}$ its universal 
covering. For $[h] = [(h_t)]\in \widetilde {Homeo_0(M,\mu)}$, and a continuous map $f : M \rightarrow \mathbb S^1,$ we lift the homotopy $fh_t - f : M \rightarrow \mathbb S^1$ to a map  $\overline{fh_t - f}$ from $M$ onto 
$ \mathbb R$. Fathi proved that the integral $\int_M\overline{fh_t - f}d\mu$ depends only on the homotopy class $[h]$ of $(h_t)$ and the homotopy class 
$\{f\}$ of $f$ in $[M,\mathbb S^1]\approx  H^1(M,\mathbb{Z})$, and that the map 
$$\widetilde{\mathfrak F}((h_t))(f):= \int_M\overline{fh_t - f}d\mu.$$
defines a homomorphism\\
$$\widetilde{\mathfrak{F}} :\widetilde {Homeo_0(M,\mu)}\rightarrow Hom(  H^1(M,\mathbb{Z}),\mathbb{R})\approx  H_1(M,\mathbb{R}). $$ This map 
induces a surjective group morphism $\mathcal{F}$ from $Homeo_0(M,\mu)$ onto a quotient of $ H_1(M,\mathbb{R})$ by a discrete subgroup.
This map is called the Fathi mass flow.

\subsubsection{Poincar\'{e} duality} Let $k$ be a fixed integer such that $0\leq k \leq 2n$. 
Consider the following bilinear map, 
$$ \langle, \rangle_P^k : H^k (M, \mathbb R)\times H^{2n -k} (M, \mathbb R)\rightarrow \mathbb R,$$
$$([\alpha],[\beta])\mapsto \int_M \alpha \wedge \beta,$$
is well defined, i.e. $\langle [\alpha],[\beta]\rangle_P^k$ does not depend on the choice of representatives 
in the cohomology class $[\alpha]$ or $[\beta]$ (this is an easy application of the Stokes
formula). Moreover this bilinear form provides an isomorphism between 
$H^k (M, \mathbb R)$ and the dual space of $H^{2n -k} (M, \mathbb R)$. 
In particular, when $\alpha \in \Omega^k(M)$ is closed and satisfies $\langle [\alpha],[\beta]\rangle_P^k = 0$ 
for all $[\beta] \in H^{2n -k} (M, \mathbb R)$  then there exists 
$\theta \in \Omega^{k-1}(M )$ such that $\alpha  = d \theta.$
\subsubsection{ Fathi's duality theorem}
It is showed in \cite{Fathi80} that the flux for volume-preserving diffeomorphisms is the Poincar\'{e} dual of Fathi's mass flow. 
Furthermore, following Fathi's \cite{Fathi80}, the latter duality result can be stated as follows. Let $\sigma$ 
denotes the canonical volume form on $\mathbb S^1$ 
given by the orientation of the circle. Then, for any function $f : M\rightarrow \mathbb S^1$, we get
$$ 
\begin{array}{ccllll}\label{EF1}
\langle Flux_\omega(\Phi), [ \dfrac{\omega^{n-1}}{(n-1)!}\wedge f^\ast\sigma] \rangle_P^1= \widetilde{\mathfrak{F}}(\Phi)(f)
\end{array}
$$
\subsection{The $C^0-$metric}\label{SC2}
 Let $Homeo(M)$ be the homeomorphisms' group of $M$ equipped with the $C^0-$ compact-open topology. This is the 
metric topology induced by the distance
  $d_0(f,h) = \max(d_{C^0}(f,h),d_{C^0}(f^{-1},h^{-1}))$
where $d_{C^0}(f,h) =\sup_{x\in M}d (h(x),f(x))$ and 
 $ d$ is a distance on $M$ induced by the Riemannian metric $g$. 
 On the space of all continuous paths $\varrho:[0,1]\rightarrow Homeo(M)$ such that $\varrho(0) = id$, 
we consider the $C^0-$topology as the metric topology induced by the metric  
$\bar{d}(\lambda,\mu) = \max_{t\in [0,1]}d_0(\lambda(t),\mu(t)).\label{d0}$

\section{Topological symplectic isotopies  \cite{BDS}, \cite{Oh-M07}, \cite{ Tchuia03}}
The introduction of topological symplectic isotopies can be motivated by the following result found in \cite{TD}. 
\begin{theorem}\label{maint}(\cite{TD}) Let $(M,\omega)$ be a closed 
symplectic manifold. 
 Let $\Phi_i = \{\phi_i^t\}$ be a sequence of symplectic isotopies,
 $\Psi = \{\psi^t\}\in Iso(M,\omega)$, and $\phi : M\rightarrow M$ be a map such that
\begin{itemize}
 \item  $(\phi_i^1)$ converges 
uniformly to $\phi$, and 
\item $l^\infty(\Phi_i^{-1}\circ\Psi)\rightarrow0,i\rightarrow\infty$.
\end{itemize}
 Then we must have $\phi = \psi^1.$
 \end{theorem}
The $L^{(1, \infty)}$ version of Theorem \ref{maint}  was recently pointed out using the positivity of the symplectic 
displacement energy (see \cite{BDS}). 
The proof of Theorem \ref{maint} given in \cite{TD} shows that the need of the positivity result of the symplectic 
displacement energy may depend on the choice of Banyaga's topology on the space of symplectic isotopies.
 The proof of Theorem \ref{maint} given in \cite{TD}
 suggests to think of the following 
questions.\\

Let $\phi$ be
a  symplectic diffeomorphism whose flux is nontrivial. Assume that $\phi$ displaces totally a given nonempty open subset $\mathcal{O}\subset M$.
Then, is 
the Hamiltonian diffeomorphism  
  arising in Hodge's decomposition of $\phi$ can totally displace the open subset in question?\\

To put Theorem \ref{maint} into further prospective, observe that it tells us to think of the following situation: if in 
Theorem \ref{maint} the sequence $\Phi_i = \{\phi_i^t\}$
 is only Cauchy in $D^2$, then
\begin{itemize}
\item What can we say about the nature of $\phi$?
 \item  Can $\phi$  be viewed as the time-one map of some continuous path?
\end{itemize}
 An attempt to answer the above questions can justify the following definitions. 
\begin{definition}\label{dfn2}(\cite{BanTchui02})
  A continuous map $\xi:[0,1]\rightarrow Homeo(M)$ with $\xi(0) = id$ is called strong symplectic isotopy if 
there exists a $D^2-$ Cauchy sequence $(F_i,\lambda_i)\subset\mathfrak{T}(M, \omega, g)$ such that 
$\bar d(\phi_{(F_i,\lambda_i)},\xi)\rightarrow0, i\rightarrow\infty$.
\end{definition}

The following result is the $L^\infty-$ context of Lemma $3.4$ found in \cite{BanTchu}. 

\begin{lemma}\label{lcontex} Let $(\phi_{(U^i, \mathcal{H}^i)}, (U^i, \mathcal{H}^i))$ be a $(D^2 + \bar d)-$Cauchy sequence. 
 Then the following holds true,
$$\max_t osc(\Delta_t(\mathcal{H}^{i} - \mathcal{H}^{i-1},\phi_{(U^{i-1},\mathcal{H}^{i-1})}))\rightarrow0, i\rightarrow\infty.$$ 
\end{lemma}
{\it Proof of Lemma  \ref{lcontex}.} The proof of Lemma \ref{lcontex} is a verbatim repetition of the proof of Lemma 3.4 
found in \cite{BanTchu}. $\Box$ 
\begin{corollary}\label{Cauchy}
 Let $\Phi_i = (\phi_i^t)$ be a sequence of symplectic isotopies which is Cauchy with respect to the metric $\bar d$ such 
that  the sequence of $1-$parameter family 
of symplectic vector fields $ X_t^i := \dfrac{d\phi_i^t}{dt}\circ(\phi_i^t)^{-1}$ is Cauchy in the norm $\|.\|^\infty$. 
Then, the sequence of $1-$parameter family 
of symplectic vector fields $Y_t^i := -(\phi_i^{-t})_\ast( X_t^i)$ is Cauchy in the norm $\|.\|^\infty$. 
\end{corollary}
{\it Proof of Corollary \ref{Cauchy}.} Put $\Phi_i = \phi_{(U^i, \mathcal{H}^i)}$. 
By definition of $Y_t^i$, we have $$\|Y^i - Y^{i + 1}\|^\infty = \max_tosc( U_t^i \circ\phi_i^t - 
U_t^{i+1} \circ\phi_{i+1}^t - \widetilde\Delta_t(\mathcal{H}^i,\Phi_i) + \widetilde\Delta_t(\mathcal{H}^{i+1},\Phi_{i+1}))$$ 
$$+ \max_t|\mathcal{H}^i - \mathcal{H}^{i+1} |$$
$$\leq \max_tosc( U_t^i \circ\phi_i^t - 
U_t^{i} \circ\phi_{i+1}^t) + \max_tosc( U_t^i \circ\phi_{i+1}^t - 
U_t^{i+1} \circ\phi_{i+1}^t)$$ 
$$ + \max_tosc(\widetilde\Delta_t(\mathcal{H}^i,\Phi_i) - \widetilde\Delta_t(\mathcal{H}^{i+1},\Phi_{i+1})) +  \max_t|\mathcal{H}^i - \mathcal{H}^{i+1} |.$$
To achieve the proof it remains to show that $\max_tosc(\widetilde\Delta_t(\mathcal{H}^i,\Phi_i) - \widetilde\Delta_t(\mathcal{H}^{i+1},\Phi_{i+1}))$ 
tends to zero when $i$ goes at infinity. To that 
end, we compute, 
$$\max_tosc(\widetilde\Delta_t(\mathcal{H}^i,\Phi_i) - \widetilde\Delta_t(\mathcal{H}^{i+1},\Phi_{i+1}))\leq 
\max_tosc(\widetilde\Delta_t(\mathcal{H}^i,\Phi_i) - \widetilde\Delta_t(\mathcal{H}^{i},\Phi_{i+1}))$$ 
$$+ \max_tosc(\widetilde\Delta_t(\mathcal{H}^i,\Phi_{1+i}) - \widetilde\Delta_t(\mathcal{H}^{i+1},\Phi_{i+1})),$$
and derive from Lemma \ref{lcontex} that 
$$\max_tosc(\widetilde\Delta_t(\mathcal{H}^i,\Phi_{1+i}) - \widetilde\Delta_t(\mathcal{H}^{i+1},\Phi_{i+1}))\rightarrow0, i\rightarrow\infty.$$ 
Lemma $3.9$ found in \cite{TD} implies that 
$$\max_tosc(\widetilde\Delta_t(\mathcal{H}^i,\Phi_i) - \widetilde\Delta_t(\mathcal{H}^{i},\Phi_{i+1}))\rightarrow0, i\rightarrow\infty.$$ This achieves 
the proof. $\Box$\\

According to Corollary \ref{Cauchy}, Definition (\ref{dfn2}) is equivalent to the following definition.
\begin{definition}\label{dfn22}
  A continuous map $\xi:[0,1]\rightarrow Homeo(M)$ with $\xi(0) = id$ is called strong symplectic isotopy if 
there exists a sequence of symplectic isotopies $ \Phi_j = (\phi_j^t)$ such that $\bar d(\Phi_j,\xi)\rightarrow0, j\rightarrow\infty$, 
and the sequence of smooth families of symplectic vector fields $ (X_t^i)$ defined by $ X_t^i := \dfrac{d\phi_i^t}{dt}\circ(\phi_i^t)^{-1}$ 
is Cauchy in $\|.\|^\infty.$
\end{definition}
We denote by $\mathcal{P}SSympeo(M, \omega) $ the space of all strong symplectic isotopies.

 It is proved in 
 \cite{BanTchui02, Tchuia03} that $\mathcal{P}SSympeo(M, \omega) $ is a group.
In particular, if the manifold is simply connected, 
then the group  $\mathcal{P}SSympeo(M, \omega) $ reduces to the group of continuous Hamiltonian flows, and the set of generators
reduces to the set of generating functions  \cite{Oh}, \cite{Viter06}. 
The set of time-one maps of all strong symplectic isotopies coincides with the group  $SSympeo(M, \omega) $, 
of all strong symplectic homeomorphisms \cite{BanTchu, Ban10c}.\\

Note that in Definition (\ref{dfn2}), if 
the sequence $(F_j,\lambda_j)$ is such that $ F_j = 0 $ for all $j$, 
then the corresponding strong symplectic isotopy is called a topological harmonic isotopy (see  \cite{Tchuia03}).\\

 For simplicity, in the rest of the present paper, 
we will sometimes say that a sequence of symplectic isotopies $(\phi_i^t)$ is Cauchy in $\|.\|^\infty$ to mean 
to mean the sequence of smooth families of symplectic vector fields $ (Z_t^i)$ defined by $ Z_t^i := \dfrac{d\phi_i^t}{dt}\circ(\phi_i^t)^{-1}$ 
is Cauchy in $\|.\|^\infty.$\\

 Some of the following questions can be found in \cite{Tchuia03}. 
\subsubsection*{Question (a)}
 Is the space of all topological harmonic flows strictly contain the space of smooth harmonic isotopies?
 \subsubsection*{Question (b)}
 What is the intersection of the set of all topological harmonic flows and  the space of all topological Hamiltonian flows? 
 
 \subsubsection*{Question (c)}
 Is any continuous path $ \gamma : t\mapsto \gamma_t $ in $SSympeo(M, \omega) $ with 
$\gamma_0 = id$, a strong symplectic isotopy?\\
 
 It is not too hard to see that the uniqueness results of Theorem \ref{URBT} and 
Theorem \ref{UGG} imply that the set of all topological harmonic flows 
 intersects the space of all 
 topological Hamiltonian flows, and the latter intersection contains a single element which is the constant path identity.\\

We start the main results of this paper by scanning the following well-known result from flux geometry  
which is due to Banyaga \cite{Ban78, Ban97} 
: ''Any symplectic 
isotopy whose flux is sufficiently small in ${H}^1(M, \mathbb{R})$ 
is homotopic relatively with  fixed extremities to a Hamiltonian isotopy''. This follows from the discreteness of the 
flux group $\Gamma$ that was proved in Ono \cite{Ono}.  
However, in the $C^0$ case, such a result is not yet known, i.e. it is 
not know any criterion that may satisfies a given ssympeotopy $\xi$ (which is not a Hameotopy) to be homotopic relatively with  fixed extremities a Hameotopy.  Analyzing the prove that was given by Banyaga \cite{Ban78, Ban97} in the 
smooth case together with the construction of Fathi's mass flow, 
  it follows from our intuition that in the $C^0$ case, a
 key ingredient of such a problem lies in Fathi's Poincar\'e duality theorem \cite{Fathi80}. To that end, we point out the following result showing 
an interplay between the symplectic dynamical 
system and the topological dynamical 
system.  
 
\begin{theorem}\label{Defor2}(Sequential-Homotopic Fathi's duality theorem) 
 Let $(M,\omega)$ be any closed symplectic manifold, and $\xi$ be a ssympeotopy whose Fathi's mass flow is trivial.  
If $(\Phi_i)$ is any sequence of symplectic isotopies which is Cauchy in $\|.\|^\infty$, and 
converges in $\bar d$ to $\xi$. Then, 
 one can extract a subsequence $(\nu_i)$ of the sequence $(\Phi_i)$ so that for each $i$, the isotopy $\nu_i$ is homotopic 
 to a Hamiltonian isotopy $\Psi_i$ (relatively with fixed extremities).   
The sequence $(\Psi_i)$ is Cauchy in $\|.\|^\infty$, and converges in $\bar d$ to a Hameotopy $\mu$.
Furthermore, the paths $\mu$  is homotopic relatively with fixed extremities to $\xi$. 
\end{theorem}

In particular, Theorem \ref{Defor2} implies that any ssympeotopies whose Fathi's mass flow 
is trivial admits its extremities in the group of Hamiltonian homeomorphisms defined in \cite{Oh}. 
Some avatars of Theorem \ref{Defor2} are underlying to the following facts.\\
 
\begin{enumerate}

\item If $\xi$ is a loop at the identity, then can we deform the sequence $(\Phi_i)$ into a sequence of loops at the identity $(\Psi_i) $ 
such that $(\Psi_i) $ is Cauchy in $\|.\|^\infty,$ and 
converges in $\bar d$ to $\xi$?
\item If $\xi$ is homotopic to the constant path identity, then can we  
extract a subsequence $(\Psi_i)$ of the sequence $(\Phi_i)$ such that for each $i$,  
the isotopy $\Psi_i$ is homotopic to the constant path identity?
\end{enumerate}

It is not too hard to see that each of the above question falls within the ambit 
of uniform approximation of (volume-preserving) homeomorphisms by
(volume-preserving) diffeomorphisms. Such problems had been studied in 
Eliashberg \cite{Elias87}, Gromov \cite{MGrom85}, 
M\"{u}ller \cite{ Mull10}, Munkres \cite{Mun1}, and Connell \cite{Conn1, Conn2}. The following result gives an affirmative answer of question $(1)$.

\begin{theorem}\label{Defor1}(Sequential Weinstein's deformation )
 Let $(M,\omega)$ be any closed symplectic manifold. Let $\xi$ be any 
strong symplectic isotopy which is a loop at the identity. 
If $(\Phi_i)$ is any sequence of symplectic isotopies which is Cauchy in $\|.\|^\infty,$ and  
converges in $\bar d$ to $\xi$.
Then, the sequence $(\Phi_i)$ 
can be deformed onto a sequence of loops $(\Psi_i)$  which is Cauchy in $\|.\|^\infty,$ and  
converges in $\bar d$ to $\xi$.
\end{theorem}

Besides of the above results, we have elaborated the following result that contributes to the comprehension of the 
topological symplectic dynamical 
systems.  In fact, it shows in particular that any strong symplectic isotopy which is a $1-$parameter group  
decomposes as composition of a smooth harmonic flow and a continuous Hamiltonian flow in the sense of Oh-M\"{u}ller.

\begin{theorem}\label{UGR}
 Let $\beta = (\beta_t)_{t\in[0,1]}$ be a strong symplectic isotopy. Assume that $\beta = (\beta_t)_{t\in[0,1]}$ is a $1-$parameter group 
( i.e: $\beta_{t + s} = \beta_t\circ \beta_s\hspace{0.2cm}\forall s,t\in [0,1]$ such that $(s + t)$ lies in $[0,1]$). 
Then, its generator $(F,\lambda)$ is time-independent. 
\end{theorem}

\subsection{Generators of ssympeotopies}
Let $\mathcal{N}^0([0,1]\times M\ ,\mathbb{R})$ be the completion of the metric space 
$\mathcal{N}([0,1]\times M\ ,\mathbb{R})$ for the $L^{\infty}$ Hofer norm, 
and $ \mathcal{PH}^1(M,g)_0$ be the completion of the metric space $\mathcal{PH}^1(M,g)$ for the uniform sup norm. Put, 
$$J^0(M, \omega, g) =:  \mathcal{N}^0([0,1]\times M\ ,\mathbb{R})\times \mathcal{PH}^1(M,g)_0,$$
and consider the following inclusion map 
$i_0: \mathfrak{T}(M, \omega, g)\rightarrow J^0(M, \omega, g).$ The map $i_0$ is uniformly continuous 
with respect of the topology induced by the metric $D^2$ on the space $\mathfrak{T}(M, \omega, g)$, 
and the natural topology of the complete metric space 
$J^0(M, \omega, g)$. Now, let $L(M, \omega, g)$ denotes the space $image(i_0),$ and 
$\mathfrak{T}(M, \omega, g)_0$ be the closure of $L(M, \omega, g)$ inside the complete metric space 
$J^0(M, \omega, g)$. That is, $\mathfrak{T}(M, \omega, g)_0$ consists of  pairs 
 $(U,\mathcal{H})$ where the mappings $(t,x)\mapsto U_t(x)$ and $t\mapsto\mathcal{H}_t$ are continuous, and for each $t$, $\mathcal{H}_t$ lies 
in  $\mathcal{H}^1(M,g)$ such that there exists a $D^2-$Cauchy sequence $(U^i,\mathcal{H}^i)\subset\mathfrak{T}(M, \omega, g)$ 
that converges  to $(U,\mathcal{H})\in J^0(M, \omega, g)$.\\ 

Observe that the sequence $(F_j,\lambda_j)$ in Definition (\ref{dfn2}) converges necessarily in the complete metric 
space $\mathfrak{T}(M, \omega, g)_0$. 
 The latter limit represents what we called in the beginning the ''generator'' of the strong symplectic isotopy $\xi$.\\
 
In the rest of this work, for short, we will often write
$(F^i,\lambda^i)\xrightarrow{L^\infty}(U,\mathcal{H})$ to mean that 
the sequence $(F^i,\lambda^i)$ converges to $(U,\mathcal{H})$ in the space $J^0(M, \omega, g)$.

 \subsection{Group structure of $GSSympeo(M, \omega, g)$, \cite{BanTchui02}}

\begin{definition}\label{Impor}(\cite{BanTchui02})
We define the set $GSSympeo(M, \omega, g)$ as the space  
of pairs $(\xi, (U,\mathcal{H}))$ where  $\xi$ is an $L^{\infty}-$ssympeotopy admitting $(U,\mathcal{H})$ as generator.
 \end{definition}
\begin{definition}
 We define the symplectic topology on the space $GSSympeo(M, \omega, g)$ as the subspace topology induced 
by its inclusion in the complete topological space $\mathcal{P}(Homeo(M),id)\times\mathfrak{T}(M, \omega, g)_0$.
\end{definition}
\subsubsection{Group structure of $GSSympeo(M, \omega, g)$}

The group 
structure on the space $GSSympeo(M, \omega, g)$ is defined as follows:
\\ For all $(\xi, (F,\lambda)), (\mu, (V,\theta))\in GSSympeo(M, \omega, g)$, their 
product is given by,
$$(\xi, (F,\lambda))\ast(\mu, (V,\theta)) = (\xi\circ\mu,(F +V\circ\xi^{-1} + \Delta^0(\theta,\xi^{-1}),\lambda + \theta)),$$ 
and the inverse of the element $(\xi, (F,\lambda))$ is given by, 
$$\overline{(\xi, (F,\lambda))} = (\xi^{-1} ,(-F\circ\xi -\Delta^0(\lambda,\xi) , -\lambda )),$$ where 
\begin{equation}\label{Q1}
 \Delta^0(\theta,\xi^{-1}) =: \lim_{L^\infty} (\widetilde\Delta(\theta^i,\phi_{(F^i,\lambda^i)}^{-1}),
\end{equation}
\begin{equation}\label{Q2}
 \Delta^0(\lambda,\xi) =: \lim_{L^\infty} (\widetilde\Delta(\lambda^i,\phi_{(F^i,\lambda^i})),
\end{equation}
 $(F^i,\lambda^i)$, and  $(V^i,\theta^i) $ are two arbitrary sequences in $\mathfrak{T}(M, \omega, g)$ such that 
$$(F^i,\lambda^i)\xrightarrow{L^\infty}(U,\mathcal{H}),$$ $$(V^i,\theta^i)\xrightarrow{L^\infty}(V,\theta),$$ 
 and 
$\widetilde\Delta_t(\lambda^i,\phi_{(F^i,\lambda^i)}^{-1})$ is the function  
$ \Delta_t(\lambda^i,\phi_{(F^i,\lambda^i)}^{-1})$ normalized.\\

\begin{remark}
 The functions defined in relations (\ref{Q1}) and (\ref{Q2}) do make sense because one derives from  
 Lemma $3.9$ found in \cite{TD} that both sequences of smooth functions $\widetilde\Delta(\lambda^i,\phi_{(F^i,\lambda^i)})$ and 
$\widetilde\Delta(\theta^i,\phi_{(F^i,\lambda^i)}^{-1})$ are Cauchy in the $L^\infty$ metric, hence converge in the 
space $\mathcal{N}^0([0,1]\times M ,\mathbb{R})$.\\

\end{remark}

For short, we will often write '' $(C^0 + L^\infty)-$topology'' to mean the $C^0-$ symplectic topology on the space $GSSympeo(M, \omega, g)$, and 
 often
$$(\lambda, (U,\mathcal{H})) = \lim_{C^0 + L^\infty} (\phi_{(U^i,\mathcal{H}^i)})$$ will mean that 
$(\lambda, (U,\mathcal{H}))$ is the limit in the space $GSSympeo(M, \omega, g)$ of a sequence 
$(\phi_{(U^i,\mathcal{H}^i)},  (U^i,\mathcal{H}^i))$ where $\phi_{(U^i,\mathcal{H}^i)}$ is the sequence of symplectic isotopies 
generated by $(U^i,\mathcal{H}^i)$. 
It is proved in \cite{BanTchui02} that $GSSympeo(M, \omega, g)$ is a topological group with respect to the $(C^0 + L^\infty)-$topology.\\ 

\begin{remark}

 Using a result found in \cite{TD} together with Fatou's lemma from measure theory, one can check that in relations (\ref{Q2}), if 
 $\xi$ is a continuous Hamiltonian flow in the sense of Oh-M\"{u}ller, then, 
$$\int_M \Delta_t^0(\theta,\xi)d\mu(x) =  0,$$ 
for each $t$, where $\mu$ is the Liouville measure. In addition, let $\xi$ be a ssympeotopy, and $\alpha$ be any closed $1-$form over $M$.  
If $(\Phi_i)$ is any sequence of symplectic isotopies which is Cauchy in $\|.\|^\infty$, and 
converges in $\bar d$ to $\xi$. Then, the mapping  following mapping is well-defined 
$$S^\alpha : \xi\longmapsto \dfrac{1}{n}\lim_{C^0 + L^\infty}\int_M\Delta_1(\alpha,\Phi_i)\omega^n.$$
In deed, If $(\Psi_i)$ is another sequence of symplectic isotopies which is Cauchy in $\|.\|^\infty$, and 
converges in $\bar d$ to $\xi$, then  
$$|\dfrac{1}{n}\int_M\Delta_1(\alpha,\Phi_i)\omega^n -\dfrac{1}{n}\int_M\Delta_1(\alpha,\Psi_i)\omega^n  | 
\leq \langle Flux_\omega(\Phi_i\circ\Psi_i^{-1}), [\alpha\wedge\omega^{n-1}] \rangle_P,$$
$$ \leq \|Flux_\omega(\Phi_i\circ\Psi_i^{-1})\|_{L^2}.\|[\alpha\wedge\omega^{n-1}]\|_{L^2},$$
for each $i$, $\langle, \rangle_P$ is the Poincar\'e scalar product, and $\|, \|_{L^2}$ stands for the $L^2-$Hodge norm on the space $H^\ast(M,\mathbb{R})$. 
Therefore, $$\lim_{C^0 + L^\infty}|\dfrac{1}{n}\int_M\Delta_1(\alpha,\Phi_i)\omega^n -\dfrac{1}{n}\int_M\Delta_1(\alpha,\Psi_i)\omega^n  | = 0, $$ since 
$\|Flux_\omega(\Phi_i\circ\Psi_i^{-1})\|_{L^2}\rightarrow0, i\rightarrow\infty.$ We get
 a linear mapping $$\mathcal{K}_0(\xi): H^1(M,\mathbb{R})\longrightarrow \mathbb{R}, [\alpha]\longmapsto S^\alpha(\xi),$$ 
for each $\xi$. Thus, $\mathcal{K}_0(\xi)$ belongs to the dual space of $H^1(M,\mathbb{R})$ which is 
$H_1(M,\mathbb{R})$. This seems to suggest that we can construct a group homomorphism 
$$ \mathcal{K}_0 : \mathcal{P}SSympeo(M, \omega)\longrightarrow H_1(M,\mathbb{R}),\xi\longmapsto \mathcal{K}_0(\xi),$$ which is similar 
to the Fathi's mass flow, and $\ker \mathcal{K}_0$ contains the group of all continuous Hamiltonian flows. 
We do not know whether $\mathcal{K}_0$ can give rise to  a $C^0-$symplectic 
invariant or not, and we do not study the above mapping in this paper. 
\end{remark}

Theorem \ref{UGG} is equivalent to the following result.

\begin{theorem}\label{UG}
 Any generator corresponds to a unique strong symplectic isotopy, i.e. 
if $(\gamma, (U,\mathcal{H})), (\xi, (U,\mathcal{H}))\in GSSympeo(M, \omega, g)$, then we must have 
$\gamma= \xi$.
\end{theorem}
\subsection{Proofs of Theorem \ref{UGG} and Theorem \ref{UGR} }
{\it Proof of Theorem \ref{UG}}. Let $(\gamma, (U,\mathcal{H}))$ and $(\xi, (U,\mathcal{H}))$ be two elements of $GSSympeo(M, \omega, g)$. 
By definition of the group $GSSympeo(M, \omega, g)$, there exist two sequences of symplectic isotopies $\phi_{(U_i,H_i)}$ and 
$\phi_{(V_i,K_i)}$ such that 
$$(\gamma, (U,\mathcal{H})) = \lim_{C^0 + L^\infty} (\phi_{(U_i,H_i)}),$$ and 
$$(\xi, (U,\mathcal{H})) = \lim_{C^0 + L^\infty} (\phi_{(V_i,K_i)}).$$
Assume that $\gamma \neq\xi,$ i.e. there exists $s_0\in ]0,1]$ such that $\gamma(s_0)\neq\xi(s_0).$ 
 Since the map $\gamma^{-1}(s_0)\circ\xi(s_0)$ belongs to $Homeo(M)$, we derive from the identity $\gamma^{-1}(s_0)\circ\xi(s_0)\neq id,$ the 
 existence  of a closed ball $B$ which is entirely moved by $\gamma^{-1}(s_0)\circ\xi(s_0)$. 
From the compactness of $B$, and the uniform convergence of the sequence $\phi_{(U_i,H_i)}^{-1}\circ\phi_{(V_i,K_i)}$ to $\gamma^{-1}\circ\xi$, we 
derive that  
\begin{equation}
 (\phi_{(U_i,H_i)}^{-s_0}\circ\phi_{(V_i,K_i)}^{s_0})(B)\cap(B) = \emptyset,\label{3E4}
\end{equation}
for all sufficiently large $i$. Relation (\ref{3E4}) implies that, 
$$0\textless e_S(B)\leq l^\infty(\phi_{(U_i,H_i)}^{-1}\circ\phi_{(V_i,K_i)}),$$ for all sufficiently large $i,$ where 
$e_S$ is the symplectic displacement energy from Banyaga-Hurtubise-Spaeth \cite{BDS}, and 
$l^\infty(.)$ represents the $L^\infty$ length functional of symplectic isotopies \cite{Ban08a}. 
On the other hand, compute  
$$ l^\infty(\phi_{(U_i,H_i)}^{-1}\circ\phi_{(V_i,K_i)}) = \max_t osc(-U_i\circ\phi_{(U_i,H_i)} + V_i\circ\phi_{(U_i,H_i)} 
 + \widetilde\Delta_t( K_i -H_i,\phi_{(U_i,H_i)}))$$
$$ + \max_t|H_i^t- K_i^t|$$
$$\leq \max_t|H_i^t- K_i^t| + \max_t osc(\widetilde\Delta_t( K_i -H_i,\phi_{(U_i,H_i)}))$$
$$ + \max_t osc(-U_i\circ\phi_{(U_i,H_i)} + U_i\circ\phi_{(U_{i+1},H_{i+1})}) $$ 
$$+ \max_t osc(-U_i\circ\phi_{(U_{i+1},H_{i+1})} + U\circ\phi_{(U_{i+1},H_{i+1})})$$
$$ + \max_t osc(-U\circ\phi_{(U_{i+1},H_{i+1})} + V_i\circ\phi_{(U_{i+1},H_{i+1})}) $$ 
$$  
+ \max_t osc(-V_i\circ\phi_{(U_{i+1},H_{i+1})} + V_i\circ\phi_{(U_i,H_i)}).$$
So, to prove that the right hand side of the above estimates tends to zero when $i$ goes at infinity, 
we only need to prove that $$\max_t osc(-U_i\circ\phi_{(U_i,H_i)} + V_i\circ\phi_{(U_i,H_i)} 
 + \widetilde\Delta_t( K_i - H_i,\phi_{(U_i,H_i)}))\rightarrow0, i\rightarrow\infty.$$  This 
follows from a straightforward application of Lemma 3.4 found in \cite{BanTchu}.  
This contradicts the positivity of the symplectic displacement energy from Banyaga-Hurtubise-Spaeth \cite{BDS}. This achieves the proof.$\square$\\

According to the above uniqueness results of strong symplectic isotopies and their 
generators vise versa, we will denote any strong symplectic isotopy $\lambda$ by $ \lambda_{(U, \mathcal{H})} $ to mean that it is 
generated by $(U,\mathcal{H})$, or equivalently 
$ \lambda_{(U, \mathcal{H})} $ is a limit
 of a sequence of symplectic isotopies 
$\phi_{(U^i,\mathcal{H}^i)}$  with respect to the $ (C^0 + L^\infty)-$topology,  i.e.  
$$\lambda_{(U, \mathcal{H})} = \lim_{C^0 + L^\infty} (\phi_{(U^i,\mathcal{H}^i)}).$$

{\it Proof of Theorem \ref{UGR}}. Let $\beta_{(F, \lambda)}$ be a strong symplectic isotopy. 
\begin{itemize}
\item Step (1). Assume that 
 $$\beta_{(F, \lambda)}^{t + s} = \beta_{(F, \lambda)}^t\circ \beta_{(F, \lambda)}^s,$$ $\forall s,t\in [0,1]$ such that $(s + t)$ lies in $[0,1]$. 
We may prove that 
$\lambda^t = \lambda^s,$ and $F^t(x) = F^s(x)$ for all $t,s\in[0,1]$, and for all $x\in M$. 
By definition of the path $t\mapsto\beta_{(F, \lambda)}^t$, we have $$\beta_{(F,\lambda)} = \lim_{C^0 + L^\infty} (\phi_{(F_i,\lambda_i)}),$$ 
where $\phi_{(F_i,\lambda_i)}$ is a sequence of symplectic isotopies. 
Observe that for each fixed $s\in[0,1]$, the sequence of symplectic maps defined by $$ \Psi^i_s (t)= \phi_{(F_i,\lambda_i)}^{(t + s)}\circ 
(\phi_{(F_i,\lambda_i)}^s)^{-1},$$ for all $t$ such that $(t +s)$ belongs to $[0,1]$, converges in $\bar d$ 
to $\beta_{(F, \lambda)}$. 
 \item  Step (2). On the other hand,  
for each $i$ compute the derivative (in $t$) of the path $ t\mapsto \Psi^i_s (t) $, and derive from the chain rule that at each time 
$t$, the tangent vector to the path $ t\mapsto \Psi^i_s (t) $ coincides with the tangent vector to the path 
$ t\mapsto \phi_{(F_i,\lambda_i)}^{(t + s)}$. That is, the isotopy $ t\mapsto \Psi^i_s (t)$ is generated by an element  $ (U^i_s, H^i_s) $ where 
 $U^i_s(t) = F_i^{t + s},$ and $ H^i_s(t) =  \lambda_i^{t + s}$ for all $t$ such that $(t +s)$ belongs to $[0,1]$, and for each $i$. 
Furthermore, the sequence of generators $  (U^i_s, H^i_s)  $  converges in the $L^\infty-$ metric to $(U_s, H_s)$ 
where $U_s(t) = F^{t + s},$ and $ H_s(t) =  \lambda^{t + s}$ for all $t$ such that $(t +s)$ belongs to $[0,1]$.
 \item  It follows from steps (1) and (2) that the element $(U_s, H_s )$ generates the strong symplectic isotopy $\beta_{(F, \lambda)}$. 
 Thus, Theorem \ref{URBT} tells us that for 
each fixed $s\in[0,1]$ we must have  
$\lambda^t = \lambda^{t + s},$ and $ F^t(x) = F^{s + t}(x)$ for all $t\in[0,1]$ such that $(t +s)$ belongs to $[0,1]$ and for all $x \in M.$
This is always true for a given $s\in[0,1]$ such that $(t +s)$ belongs to $[0,1]$, i.e. we have  
$\lambda^t = \lambda^0,$ and $F^t(x) = F^0(x)$ for all $t\in[0,1]$, and for all $x\in M.$ This achieves the proof. $\Box$

\end{itemize}

 As we said in the beginning, Theorem \ref{UGR} suggests that any strong symplectic isotopy which is a $1-$parameter 
group decomposes into the composition of smooth harmonic 
flow and a 
continuous Hamiltonian flow in the sense of Oh-M\"{u}ller \cite{Oh-M07}.

\subsubsection*{Question (e)} Is $1-$parameter group any strong symplectic isotopy which decomposes into the composition of smooth harmonic 
flow and a 
continuous Hamiltonian flow?\\

We have the following fact.

\begin{lemma}\label{UGR-1} Let $
\lambda_{(U, \mathcal{H})}^t$, $t\in [0,1]$ be any strong symplectic isotopy. For each fixed $s\in [0,1)$, the path 
$ t\mapsto \lambda_{(V, \mathcal{K})}^t := \lambda_{(U, \mathcal{H})}^{(t + s)}\circ  (\lambda_{(U, \mathcal{H})}^s)^{- 1} $ 
is a strong symplectic isotopy generated by 
$(V, \mathcal{K})$ where $ V(t,x) = U(t + s, x),$ and $ \mathcal{K}_t= \mathcal{H}_{(t + s)}$ for all $t\in [0,1 -s]$, and for all $x\in M$.
\end{lemma}
{\it Proof of Lemma \ref{UGR-1}}. Assume that $\lambda_{(U, \mathcal{H})} = \lim_{C^0 + L^\infty} (\phi_{(U^i,\mathcal{H}^i)}).$ 
For each fixed $s\in [0,1)$, consider the sequence $ \phi_{(V^i, \mathcal{K}^i)}$ of symplectic isotopies defined by   
$$ \phi_{(V^i, \mathcal{K}^i)}^t = \phi_{(U^i, \mathcal{H}^i)}^{(t + s)}\circ  (\phi_{(U^i, \mathcal{H}^i)}^s)^{- 1},$$
for all $t$, and for each $i$.  
Compute the derivative (in $t$) of the path $ t\mapsto \phi_{(V^i, \mathcal{K}^i)}^t $, and derive from the chain rule that at each time $t$, 
the tangent vector to the path $ t\mapsto \phi_{(V^i, \mathcal{K}^i)}^t $ coincides with the tangent vector to the path
 $ t\mapsto \phi_{(U^i, \mathcal{H}^i)}^{(t + s)}$, or equivalently we get $ V^i(t) = U^i( t+s) ,$ and 
$ \mathcal{K}^i_t =  \mathcal{H}^i_{(t + s)}$ for all $t\in[0,1]$ such that $(t +s)$ belongs to $[0,1]$. 
A straightforward computation implies that the sequence of symplectic 
isotopies $ t\mapsto \phi_{(V^i, \mathcal{K}^i)}^t  $ converges in $\bar d$ to 
$\lambda_{(U, \mathcal{H})}^{(t + s)}\circ  (\lambda_{(U, \mathcal{H})}^s)^{- 1} $, as well as the sequence 
of generators 
$ (V^i,\mathcal{K}^i) $ converges in the $L^\infty$ topology to an element $ (V, \mathcal{K}) $ such that 
$ V(t,x) = U(t+ s, x) ,$ and $ \mathcal{K}_t= \mathcal{H}_{(t + s)}$ for all $t\in [0,1 -s]$, and for all $x\in M$. 
That is, $\lambda_{(V, \mathcal{K})}$ is a strong symplectic isotopy generated by $(V, \mathcal{K}),$ i.e.  
$\lambda_{(V, \mathcal{K})}= \lim_{C^0 + L^\infty} (\phi_{(V^i,\mathcal{K}^i)}).$ This completes the proof. $\Box$\\

In the following, by $\sharp$ we denote the natural isomorphism induced by the 
symplectic form $\omega$ from cotangent bundle $TM^\ast$ to tangent bundle $TM$.\\

We will need the following lemma.\\
\begin{lemma}\label{LDefor2}(Sequential deformation)
 Let $(M,\omega)$ be any closed symplectic manifold. Let $(Z^i_t)$ be a $ \| .\|^\infty-$Cauchy sequence of harmonic 
vector fields. For each $i$, set 
$$Z^{(s,t)}_i = tZ_{st}^i - 2s(\int_0^t(\iota(Z_u^i)\omega)du)^\sharp,$$
$$ Y_i^t = -\int_0^tZ^i_u du.$$
We have the following properties.
\begin{enumerate}

\item For each $i$, $(Y_i^t)_t$ is a smooth family of harmonic vector fields, 
the sequence $(Y_i^t)_t$ is Cauchy in $\|. \|^\infty$, and the  sequence smooth paths 
generated by $(Y_i^t)_t$ converges in $\bar d$. The latter limit is a strong symplectic isotopy.
\item For each fixed $t$, let $(\theta_{s,t}^i)_s$ denotes the flow generated by $ Y_i^t$. 
Then, $(\theta_{s,t}^i)_s$ converges in $\bar d$. The latter limit is a strong symplectic isotopy.
\item For each $i$, for each fixed $t$, 
sequence of family of symplectic vector fields $(Z^{(s,t)}_i)_s$ converges in $\|. \|^\infty$, and the sequence of smooth paths
 $(G_{(s,t)}^i)_s$ generated by  $(Z^{(s,t)}_i)_s$ converges in $\bar d$. The latter limit is then obviously a strong symplectic isotopy.
\item For each fixed $s$, the 
sequence of family of symplectic vector fields defined by $V_{(s,t)}^i = \dfrac{d}{dt}G_{(s,t)}^i((G_{(s,t)}^i)^{-1})$ 
is Cauchy in $\|. \|^\infty$.
\end{enumerate}
\end{lemma}
\

{\it Proof of Lemma \ref{LDefor2}}.  This proof 
is subtle, in view of this fact we shall proceed step by step.
 Let $(Z^i_t)$ be a sequence of harmonic  vector fields which 
 is Cauchy in $\| .\|^\infty.$  
\begin{itemize}
\item (1). For each $i$, for each $u\in[0,t]$, $Z^i_u$ is harmonic, and 

the map $t\mapsto -\int_0^tZ^i_u du$ is smooth, i.e. $Y^t_i$ is harmonic. 
Compute, $$\|(Y_i^t)_t - (Y_{1 +i}^t)_t \|^\infty\leq \|(Z^i_t)_t - (Z^{1 +i}_t)_t \|^\infty,$$
for each $i$, and conclude that the right hand side tends to zero when $i$ goes at infinity. Since 
the sequence $(Y_i^t)_t$ is Cauchy in $\|. \|^\infty$, we derive that the sequence of smooth 
family of harmonic $1-$forms  $(\mathcal{H}_i^t)_t$ defined by $\iota(Y_i^t)\omega = \mathcal{H}_i^t$ is 
Cauchy in $\mathcal{PH}^1(M,g)$ with respect to the norm $\|. \|^\infty$. Hence,  the latter converges in the complete metric space 
$\mathcal{PH}^1(M,g)_0$ to a continuous family of smooth harmonic $1-$forms  $(\mathcal{H}^t)$. 
Since $(\mathcal{H}^t)^\sharp$ is a continuous family of
 vector fields, we derive from \cite{Tchuia03} (Lemma 5.1, \cite{Tchuia03}) or \cite{Ral-Rob} that 
the sequence of smooth paths 
generated by $(Y_i^t)_t$ converges in $\bar d$ to a continuous family of smooth diffeomorphisms. 
Therefore, the celebrated rigidity theorem from Eliashberg-Gromov tells us that the 
latter limit is a continuous family of smooth symplectic diffeomorphisms. 
\item (2). For each fixed $t$, since the sequence of harmonic vector fields 
 $ (Y_i^t)_i$ is Cauchy in $\|. \|^\infty$, we derive as in item (1) that 
 the sequence $(\theta_{s,t}^i)_s$ of flows generated by the sequence 
$ (Y_i^t)_i$ converge in $\bar d$. 
 \item (3). It is not too hard to derive from the assumption 
that for each fixed $t$, the sequence symplectic vector fields $(Z^{(s,t)}_i)_s $ is Cauchy in $\|. \|^\infty$  
since for each fixed $t$, a straightforward calculation leads to the following estimate 
$$\| (Z^{(s,t)}_{i+ 1})_s  - (Z^{(s,t)}_i)_s\|^\infty\leq 3\|(Z_{u}^i) - (Z_{u}^{i+1})\|^\infty.$$ 
On the other hand,since the sequence of harmonic vectors fields $(Z^{(s,t)}_i)_s $ is Cauchy in $\|. \|^\infty$, 
one uses the same arguments as in item (1) to derive that 
the sequence of symplectic isotopies $(G_{(s,t)}^i)_s$ generated by $(Z^{(s,t)}_i)_s $ converges in $\bar d$ to a continuous family 
$(G_{(t,s)})_s$ of 
symplectic diffeomorphisms. That is, for each fixed $t$, the continuous family $(G_{(t,s)})_s$ of 
symplectic diffeomorphisms is a strong symplectic isotopy so that the map 
 $(s,t)\mapsto G_{(s,t)}$ is continuous. 
It follows from the above that for each fixed $s$, the sequence of symplectic isotopies $(G_{(s,t)}^i)_t$ converges 
uniformly to the continuous family $(G_{(t,s)})_t$ of 
symplectic diffeomorphisms.

\item (4). The main argument used in the following can be found in \cite{Ban08a}.
 Let $symp(M,\omega)$ be the space of symplectic vector fields.
 Consider $\mathfrak{X}$ to be the space of all smooth curves, $$c:I=[ 0,1]\rightarrow symp(M,\omega),$$ 
with $c(0) = 0$ endowed with the norm
$\|.\|^\infty.$ Then, we equip the product space $\mathfrak{X}\times I$ with the following distance
$$\delta((c,s),(c',s')) = ((\|c - c'\|^\infty)^2 + |s - s'|^2)^{1/2}.$$ On the other hand, consider 
$\mathfrak{N}$ to be the space of all smooth functions $u:I\times I \rightarrow symp(M,\omega)$ endowed with the norm 
$$\|u\|_0 = \sup_{s,t}\|u(s,t)\|.$$ We have the following smooth mappings, 
$$
\begin{array}{lllc}
 A_s : c(t)\mapsto tc(st) - 2s(\int_0^t(i_{c(u)}\omega) du)^\sharp = U_{s,t} \\
I_s :U_{s,t}\mapsto G_{s,t} \\
\partial_t : G_{s,t}\mapsto \frac{\partial}{\partial t} G_{s,t},
\end{array}
$$
which induce the following Lipschitz map $\mathcal{R} :\mathfrak{X}\times I\rightarrow \mathfrak{N}$ where $\mathcal{R} = \partial_t\circ I_s\circ A_s.$
Observe that for each $i$, the 2-parameter family of vector fields $V_{s,t}^i$ is the image of the couple 
 $ ((Z_{t}^i), s)$ by the Lipschitz continuous mapping $\mathcal{R}$.  
Then, it follows from the Lipschitz uniform continuity that for all fixed $s$, the sequence of symplectic vector fields $(V_{s,t}^i)_t$ 
is Cauchy in $\|. \|^\infty$ because  
\begin{eqnarray}
\sup _t| V_{s,t}^i - V_{s,t}^{1 + i}| & =& \delta (\mathcal{R}((Z_{t}^i), s), \mathcal{R}((Z_t^{1 +i}),s))\cr\cr 
&\leq& \kappa \|(Z_{t}^i) -(Z_t^{1 +i})\|^\infty \cr\cr
\end{eqnarray}
where $\kappa$ is the Lipschitz constant of the map $\mathcal{R}$. The right hand side of the above estimate tends 
to zero when $i$ goes at infinity. This achieves the proof. $\Box$ 

\end{itemize}

\section{An enlargement of the first Calabi's invariant}

Note that in view of the uniqueness results of Theorem \ref{URBT} and 
 Theorem \ref{UGG} the groups $\mathcal{P}SSympeo(M, \omega)$ and $ GSSympeo(M, \omega, g)$ are isomorphic, i.e. 

$$ GSSympeo(M, \omega, g) \approx\mathcal{P}SSympeo(M, \omega).$$
So, we identify $GSSympeo(M, \omega, g)$ with $\mathcal{P}SSympeo(M, \omega)$, and 
refer to the topology on $\mathcal{P}SSympeo(M, \omega)$ as the $(C^0 + L^\infty)-$topology. Under this identification, 
$\mathcal{P}SSympeo(M, \omega)$ can be viewed as a topological group.
\begin{definition} We define 
 the symplectic topology on the space $SSympeo(M, \omega)$ to be the strongest topology 
which makes the mapping 
 $$ ev : \mathcal{P}SSympeo(M, \omega)\rightarrow  SSympeo(M, \omega),$$
 $$ \xi_{(F,\lambda)}\mapsto \xi_{(F,\lambda)}^1,$$
becomes continuous with respect to the $(C^0 + L^\infty)-$topology on $\mathcal{P}SSympeo(M, \omega)$.
\end{definition}
Denote by $\mathfrak{SS}ympeo(M,\omega)$ the group $SSympeo(M, \omega) $
  equipped 
with the symplectic topology. By definition, the map  $$ev : \mathcal{P}SSympeo(M, \omega)\rightarrow \mathfrak{SS}ympeo(M,\omega)$$ 
is surjective, continuous and open. In fact, the openness of $ev$ is proved as follows. 
 Pick an open subset $\mathcal{O}$ inside the space $\mathcal{P}SSympeo(M,\omega)$. We have to prove that $ev(\mathcal{O})$ is an open 
subset inside the topological space $\mathfrak{SS}ympeo(M,\omega)$. By definition of the topological structure on the space $\mathfrak{SS}ympeo(M,\omega)$, 
it suffices to prove that the subset $ev^{-1}(ev(\mathcal{O}))$ is open in $\mathcal{P}SSympeo(M,\omega)$.\\ Let 
$\gamma$ be a loop inside the topological group $\mathcal{P}SSympeo(M,\omega)$, and let $L_\gamma$ denotes the 
left translation by $\gamma$ inside 
the topological group $\mathcal{P}SSympeo(M,\omega)$. We have $ev(L_\gamma(\mathcal{O})) = ev(\mathcal{O})$ which 
implies that $L_\gamma(\mathcal{O})\subset ev^{-1}(ev(\mathcal{O}))$. Then, by taking the union over all the 
loops $\gamma$ in $\mathcal{P}SSympeo(M,\omega)$ we obtain\\ 
$\bigcup_{\gamma}L_\gamma(\mathcal{O})\subset ev^{-1}(ev(\mathcal{O})).$ 
Since translations are open in any topological group, it follows that  
 $\bigcup_{\gamma}L_\gamma(\mathcal{O})$ is open as an arbitrarily 
collection of open sets (we refer to \cite{TPG-66}, \cite{TPG-77} for further comprehension of topological groups). On the another hand, let $\alpha$ be 
an element that belongs to  $ev^{-1}(ev(\mathcal{O}))$. By characterization of any element therein $ev^{-1}(ev(\mathcal{O}))$, 
there exists an element $\beta$ in $\mathcal{O}$ such that $\alpha = ev^{-1}(ev(\beta))$. That is, $ev(\alpha) = ev(\beta)$ since $ev$ is surjective. 
From the identity $\alpha = (\alpha\circ\beta^{-1})\circ\beta$ we derive that 
$\alpha \in L_{\alpha\circ\beta^{-1}}(\mathcal{O})$ since $\alpha\circ\beta^{-1}$ is a loop at the identity. That is, 
$ev^{-1}(ev(\mathcal{O}))\subset\bigcup_{\gamma}L_\gamma(\mathcal{O}).$ Finally, 
$$ev^{-1}(ev(\mathcal{O})) = \bigcup_{\gamma}L_\gamma(\mathcal{O}),$$ is open. As a consequence of the above fact we see 
that $\mathfrak{SS}ympeo(M,\omega)$ is a topological group. 
\subsection{The mass flow for strong symplectic isotopies}
In the following subsection, we compute the mass flow for 
strong symplectic isotopies.\\

Let $(\Phi_i)$ be a sequence of symplectic isotopies which converges in $\bar d$ to $\lambda$. According to Fathi's, the assignation 
$\lambda \mapsto \widetilde{\mathfrak{F}}(\lambda)$ is continuous with the respect to the uniform topology, i.e.
 $$\lim_{\bar d} (\widetilde{\mathfrak{F}}(\Phi_i)(f)) = \widetilde{\mathfrak{F}}(\lambda)(f),$$ 
for all continuous mapping $f : M \rightarrow \mathbb S^1$. This 
implies that $$\lim_{\bar d}\langle Flux_\omega(\Phi_i), [ \dfrac{\omega^{n-1}}{(n-1)!}\wedge f^\ast\sigma]\rangle_P^1 
= \widetilde{\mathfrak{F}}(\lambda)(f),$$
 for all continuous mapping $f : M \rightarrow \mathbb S^1$ (see Section 3.3 of the present paper). 
But, in the latter equality one cannot permute the limit in $\bar d$ with the integral since the convergence  of 
the sequence $(\Phi_i)$ in $\bar d$ does not guarantee
 that the sequence of cohomological classes $(Flux_\omega(\Phi_i))$  converges in $ H^1(M,\mathbb{R})$ with respect to the topology induced by the vector 
space structure. For instance, 
 consider the following continuous mappings, 
$$P_{2} : \mathcal{P}SSympeo(M, \omega)\rightarrow \mathcal{PH}^1(M,g)_0,$$ 
$$\xi_{(F,\theta)}\mapsto \theta,$$ 
  
$$Q_0 : \mathcal{P}SSympeo(M, \omega)\rightarrow \mathcal{P}(Homeo(M),id),$$ 
$$\xi_{(F,\theta)}\mapsto \xi_{(F,\theta)}.$$ 
Since any symplectic isotopy $\phi_{(U,\mathcal{H})} $ belongs to $\mathcal{P}SSympeo(M, \omega)$, 
 we derive that, 
\begin{equation}\label{EF2}
 \langle[\int_0^1 P_{2}(\phi_{(U,\mathcal{H})})(t)dt], 
[ \dfrac{\omega^{n-1}}{(n-1)!}\wedge f^\ast\sigma]\rangle_P^1 = 
\widetilde{\mathfrak{F}}(Q_0(\phi_{(U,\mathcal{H})}))(f),
\end{equation}
where $$[\int_0^1 P_{2}(\phi_{(U,\mathcal{H})})(t)dt] =  \int_0^1 [\mathcal{H}_t]dt = Flux_\omega(\phi_{(U,\mathcal{H})}).$$
Let $\lambda_{(U,\mathcal{H})}$ be a strong symplectic isotopy such that 
$$\lambda_{(U,\mathcal{H})} = \lim_{(C^0 + L^\infty)}(\phi_{(U^i,\mathcal{H}^i)}).$$ 
 According to equation (\ref{EF2}), we have,
$$ \langle[\int_0^1 P_{2}(\phi_{(U^i,\mathcal{H}^i)})(t)dt], 
[ \dfrac{\omega^{n-1}}{(n-1)!}\wedge f^\ast\sigma]\rangle_P^1 = 
\widetilde{\mathfrak{F}}(\phi_{(U^i,\mathcal{H}^i)})(f),$$
for each $i$, and taking the limit with respect to the $(C^0 +L^\infty)-$topology in both sides leads to,
$$
\begin{array}{lllccc}
 \widetilde{\mathfrak{F}}(\lambda_{(U,\mathcal{H})})(f) 
&=&  \lim_{(C^0 + L^\infty)}\langle[\int_0^1 P_{2}(\phi_{(U^i,\mathcal{H}^i)})(t) dt], 
[  \dfrac{\omega^{n-1}}{(n-1)!}\wedge f^\ast\sigma]\rangle_P^1\cr\cr
&=&   \langle \lim_{L^\infty}Flux_\omega((\phi_{(U^i,\mathcal{H}^i)}), [  \dfrac{\omega^{n-1}}{(n-1)!}\wedge f^\ast\sigma]\rangle_P^1\cr\cr
&=&   \langle[ \int_0^1 \mathcal{H}_tdt], [  \dfrac{\omega^{n-1}}{(n-1)!}\wedge f^\ast\sigma]\rangle_P^1. 
\end{array}
$$
Therefore, it follows from the above estimates that the mass flow of any strong symplectic isotopy $ \lambda_{(U,\mathcal{H})}$ 
is given by : 
\begin{equation}\label{FmF}
\widetilde{\mathfrak{F}}(\lambda_{(U,\mathcal{H})})(f) =  
\dfrac{1}{(n-1)!}\int_M  ( \int_0^1 \mathcal{H}_tdt)\wedge \omega^{n-1}\wedge f^\ast\sigma,
\end{equation}
for all continuous mapping $f : M \rightarrow \mathbb S^1$.  
In particular, since any continuous Hamiltonian flow $\xi$ (in the sense of Oh-M\"{u}ller) can be written as $\xi = \xi_{(U, 0)}$, 
we derive from equation (\ref{FmF}) that the mass flow of any continuous Hamiltonian flow is trivial. This 
agrees with a result that was prove by Oh-M\"{u}ller \cite{Oh-M07} 
asserting that the mass flow of any continuous Hamiltonian flow is trivial. 
We are now ready to prove  Theorem \ref{Defor2}. 

\subsection{Proofs of Theorem \ref{Defor2} and Theorem \ref{Defor1}}
Recall that Theorem \ref{Defor2} 
states that if $\xi$ is a strong symplectic isotopy such that $ \widetilde{\mathfrak{F}}(\xi)(f) = 0,$ 
for any function $f : M\rightarrow \mathbb S^1$, 
then $\xi$ is homotopic relatively with fixed extremities to a continuous Hamiltonian flow.\\

{\it Proof of Theorem \ref{Defor2}}. 
Let $\xi$ be a continuous symplectic 
flow whose Fathi's mass flow is trivial, and $(\Phi_i)$ be a sequence of symplectic isotopies 
such that $\xi = \lim_{C^0 + L^\infty}(\Phi_i).$ For each $i$, 
set $\Phi_i = (\phi_t^i),$ and 
let $\phi_t^i = \pi_t^i\circ\alpha_t^i$ be the Hodge decomposition of the isotopy $(\phi_t^i)$ 
where $(\pi_t^i)$  is a harmonic isotopy, and  $(\alpha_t^i)$ is a Hamiltonian isotopy. In 
view of Hodge's decomposition 
of strong symplectic isotopies, one can write the ssympeotopy $\xi$ as follows :$$\xi(t) = \pi_t\circ\alpha_t,$$ where $(\pi_t)$ 
is a continuous harmonic flow, and $(\alpha_t)$ is 
continuous Hamiltonian flow in the sense of Oh-M\"{u}ller (see \cite{Tchuia03}, \cite{Oh-M07}). 
\begin{itemize}
 \item Step (a). By assumption, the 
 Fathi's mass flow of $\pi$ is trivial. Thus, the continuity of 
the mapping
$\lambda \mapsto \widetilde{\mathfrak{F}}(\lambda)$ tells us that 
there exists a large integer $i_0$ such that for all $ i\textgreater i_0$, the 
mass flow $ \widetilde{\mathfrak{F}}(\rho_i)$ is arbitrarily small. 
Next, we derive from Fathi's Poincar\'{e} duality 
theorem that one can make the flux of $\rho_i= (\rho_t^i)$ arbitrarily small. Assume this done. 
Then, under this assumption, we derive from Theorem \ref{Fgeo} that 
for all $i \textgreater i_0$, the isotopy $\rho_i = (\rho_t^i)$  is homotopic relatively with fixed endpoints 
 to a Hamiltonian isotopy $(\varphi^i_t)$ constructed as in the proof of Theorem \ref{Fgeo}. More precisely, 
it follows from the proof of Theorem \ref{Fgeo} that for each $t$, we have  $$\varphi^i_t = \theta_{1,t}^i\circ\rho_t^i,$$
for each $i$, 
 where for each fixed $t$,
$(\theta_{s,t}^i)_s$ is the flow generated by the symplectic 
vector field $Y_t = - (\int_0^t\iota(\dot{\rho}_u^i )du)^\sharp.$ 
One derives from Lemma \ref{LDefor2} that 
the sequence $(\varphi^i_t)$ converges to a continuous Hamiltonian flow $(\varphi_t)$. 
\item Step (b). Consider the sequence $(\nu_i) = (\Phi_i)_{i \textgreater i_0}$, and derive from step (a) that for 
each $i$ the path $\nu_i$ is homotopic to relatively with fixed extremities to $\varphi^i_t\circ \alpha_t^i$ where 
the isotopy  $(\alpha_t^i)$ is the Hamiltonian part in Hodge's decomposition of 
the isotopy $\Phi_i = (\phi_t^i)$. As a first glance, observe that 
the sequence $(\varphi^i_t\circ \alpha^i)_t$ 
converges  in the $(C^0 + L^\infty)-$topology to the continuous Hamiltonian flow $ \mu : t\mapsto \varphi_t\circ \alpha_t$ where  
 the isotopy  $(\alpha_t)$ is the Hamiltonian part in Hodge's decomposition of 
the ssympeotopy $\xi$ as indicated in the beginning of this proof. Next, compute 
$$\mu(1) = \lim_{C^0}( \varphi^i_1\circ \alpha_1^i) = \lim_{C^0}( \rho_1^i\circ \alpha_1^i) = \lim_{C^0}( \phi_i^1 ) = \xi(1).$$
\item Step (c).  Lemma \ref{LDefor2}-(1), (2) suggests that
\begin{enumerate}
 \item for each fixed $t$, the sequence of symplectic isotopies $(\theta_{s,t}^i)_s$ converges in $\mathcal{P}SSympeo(M, \omega)$ to a 
strong symplectic isotopy $ (\theta_{(s,t)})_s$,
\item for each fixed $s$, the sequence of symplectic isotopies $(\theta_{s,t}^i)_t$ converges in $\mathcal{P}SSympeo(M, \omega)$ to a 
strong symplectic isotopy $ (\theta_{(s,t)})_t$.
\end{enumerate}
It follows from the above facts that for all $(s,t)\in [0,1]\times[0,1]$, the $C^0$ limit of the sequence of symplectic diffeomorphisms 
$(\theta_{s,t}^i)$ exists, and the latter limit lies in $SSympeo(M, \omega)$. For all $(s,t)\in [0,1]\times[0,1]$, 
set $$\lim _{C^0}(\theta_{s,t}^i) = \theta_{(s,t)}.$$ 
By construction, the two paths $(\varphi_t)$ and  $t\mapsto \theta_{(1,t)}\circ\pi_t$ coincide, and 
 $\theta_{(s,0)} = \theta_{(s,1)} = id_M = \theta_{(0,t)}.$ Thus, 
$$ \mu : t\mapsto \varphi_t\circ \alpha_t = \theta_{(1,t)}\circ\pi_t\circ \alpha_t ,$$ is a continuous 
Hamiltonian flow.
Finally, we define a homotopy  between $\xi$ and $\mu$ as follows. 
$$ H : [0, 1]\times [0, 1]\rightarrow \mathfrak{SS}ympeo(M,\omega),$$ 
$$ (s,t)\mapsto \theta_{(s,t)}\circ\pi_t\circ\alpha_t,$$
such that
$H(s, 0) = id$,  $H(s, 1) = \xi(1) = \mu(1)$, $H(0, t) = \xi(t)$  and $H(1, t) = \mu (t)$.
This completes the proof. $\Box$\\
\end{itemize}

{\it Proof of Theorem \ref{Defor1}}. Since the sequence $\Phi_i$ converges in $(C^0 + L^\infty)-$topology to the loop $\xi$, 
we derive that the sequence of time one maps 
$\Phi_i(1)$
converges uniformly to the constant map identity. In view of Hodge's decomposition theorem 
of symplectic isotopies, it is known that for each $i$, the diffeomorphism $\Phi_i(1)$ decomposes 
as $\Phi_i(1) = \rho_i\circ\varrho_i$ where $\rho_i$ is a harmonic diffeomorphism and $\varrho_i$ 
a 
Hamiltonian diffeomorphism. Hence, the bi-invariance of the metric $\bar d$ suggests that 
$$ d_{C^0}(\rho_i^{-1}(1), \varrho_i(1)) =  d_{C^0}(\Phi_i(1), id) \rightarrow0, i\rightarrow\infty,$$  i.e. the 
sequences $(\rho_i^{-1}(1))$ and $(\varrho_i(1))$ of time-one maps  converges uniformly to the same limit $\rho$. 
But, a result found in  \cite{Tchuia03} (Lemma 5.1, \cite{Tchuia03}) shows that 
the sequence $(\rho_i)$ always converges in $\bar d$ to a continuous path $(\rho_t)$ in $Symp_0(M,\omega)$. 
Thus, it follows from the above statements 
that the sequence $(\varrho_i(1))$ converges uniformly to the symplectic diffeomorphism $\rho^{-1}_1$. 
That is, the sequence $(\Phi_i(1))$ converges uniformly to the 
diffeomorphism $\rho\circ\rho^{-1} = id$. Hence, we derivative from 
the celebrated rigidity result dues to Eliashberg \cite{Elias87} that 
we can find a large integer $i_0$ so that for $\forall i\textgreater i_0$,
the diffeomorphism $\Phi_i(1)$ lies in a small   
$C^\infty$ neighborhood $ \mathcal{W}(i,id)$ of the constant map identity in $Symp_0(M,\omega)$ with the following property:
for 
all $i,j\textgreater i_0$ such that $j\textgreater i$, we have $$ \mathcal{W}(j,id)\subset\mathcal{W}(i,id).$$
But, in view of a result that was proved by Weinstein \cite{AWein}, the group $Symp_0(M,\omega)$ is 
locally contractible, and then locally connected by smooth arcs with respect to the $C^\infty$ compact-open topology. 
So, for all $i\textgreater i_0$, each $\Phi_i(1)$ can be connected 
to the identity through a smooth symplectic isotopy $\beta_i$ with the following properties : 
\begin{itemize}
\item for each $t$, $ \beta_i^t\in \mathcal{W}(i, id),$ 
\item the $L^\infty$ Banyaga's Hofer-like 
length of the isotopy $\beta_i$  is less or equal to the $L^\infty$ Banyaga's Hofer-like 
length of any other symplectic isotopy that connects $\beta_i^1= \Phi_i(1)$ to the identity.

\end{itemize}
 The 
latter process generates automatically a sequence $(\beta_i)$ of small symplectic isotopies 
which Cauchy in $\|.\|^\infty$ since $l^\infty(\beta_i) \rightarrow0, i\rightarrow\infty$. 
 Now, set 
$\Psi_i = \Phi_i\circ \beta_i^{-1}$ for all $i\textgreater i_0$. This defines 
a sequence loops at the identity which is Cauchy in $\|.\|^\infty$. The latter sequence 
converges 
 in $\bar d$ to $\xi.$ This achieves the proof.$\Box$

\subsection{Strong Calabi's invariant}
The formula of mass flow for strong symplectic isotopies suggests that we can define a $C^0-$ flux geometry 
for strong symplectic homeomorphisms so that Fathi's Poincar\'e duality theorem continues to hold (at least when $M$ is of type Lefschetz). Notice 
that $M$ is of type Lefschetz, if  
the mapping $$\cup^\omega : H^{1}(M,\mathbb{R})\rightarrow H^{2n-1}(M,\mathbb{R}), [\alpha]\longmapsto  [\alpha]\wedge[\omega^{n-1}],$$
is an isomorphism. 
 For this purpose, consider the following mapping,
$$ \widetilde{Cal_0} : \mathcal{P}SSympeo(M, \omega) \xrightarrow{}   H^{1}(M,\mathbb{R}) ,$$
$$ \xi_{(F,\theta)}\mapsto \int_0^1[\theta_t]dt.$$
The map $\widetilde{Cal_0}$ is a well-defined, continuous and surjective
 because of the one-to-one correspondence 
between the space of ssympeotopies and that of their generators. Furthermore, the way that 
the group structure on  $GSSympeo(M, \omega, g) $ is defined suggests that the map $\widetilde{Cal_0}$ is a  group homomorphism. 
In the rest of this work, we will sometimes use the terminology  ''strong flux'' to mean the mapping $\widetilde{Cal_0}$.
\subsubsection{Homotopic invariance of strong flux }
We equip the group of all strong symplectic isotopies with the 
  equivalent relation $\sim $ defined as follows:

Let $ \xi_{(F_i,\theta_i)} $, $ i = 1, 2 $ 
be two strong symplectic isotopies. We 
shall write $$\xi_{(F_2,\theta_2)}\sim \xi_ {(F_2,\theta_2)},$$ 
if and only if $\xi_{(F_1,\theta_1)}$ is homotopic relatively with 
 fixed extremities to $\xi_{(F_2,\theta_2)}$.\\

 We shall write  $\{\xi_{(F,\theta)}\}$ to represent the equivalent class of the ssympeotopy 
$ \xi_{(F,\theta)}$ with respect to the equivalent relation $\sim $. 
Let $\widetilde{\mathcal{P}SSympeo(M, \omega)}$ denotes the quotient space of $\mathcal{P}SSympeo(M, \omega)$ with respect to the equivalent 
relation $ \sim $. If  $ SSympeo(M, \omega)$ is locally path connected, then $\widetilde{\mathcal{P}SSympeo(M, \omega)}$ is the universal cover 
of $ SSympeo(M, \omega)$.

\begin{proposition}\label{P17} The map $ \widetilde{Cal_0 }$ does not depend on the choice of a 
representative in the equivalent class of any strong symplectic isotopy, provided $M$ is Lefschetz.\\

\end{proposition}
{\it Proof of Proposition \ref{P17}}. 
  Let  $ \xi_{(F_i,\theta_i)} $, $ i = 1, 2 $ be two strong symplectic isotopies such that  
  $\xi_{(F_2,\theta_2)}\sim \xi_ {(F_2,\theta_2)}$. 
set $$ \xi_ {(F_2,\theta_2)}= \lim_{C^0 + L^\infty}(\phi_ {(F_2^i,\theta_2^i)}),$$
and $$ \xi_ {(F_1,\theta_1)}= \lim_{C^0 + L^\infty}(\phi_ {(F_1^i,\theta_1^i)}).$$
This is equivalent to say that  
$ \xi_{(F_2,\theta_2)}\circ  \xi_{(F_1,\theta_1)}^{-1}$ is homotopic relatively with fixed endpoints to the constant 
path identity. Therefore, the Fathi's mass flow of the loop $ \xi_{(F_2,\theta_2)}\circ  \xi_{(F_1,\theta_1)}^{-1}$  is trivial. But it is not difficult 
to see that the mass flow of $\xi_{(F_2,\theta_2)}\circ  \xi_{(F_1,\theta_1)}^{-1}$ reduce to that of $ \xi_{(0,\theta_2 - \theta_1)}$. That is, 
the mass flow of  $ \xi_{(0,\theta_2 - \theta_1)}$ is trivial. Since the Fathi's mass flow is continuous with respect to the $C^0$ topology, we derive that 
for all $i$ sufficiently large, the mass flow of the harmonic isotopy $ \phi_ {(0,\theta_2^i - \theta_1^i)}$ is arbitrarily small since 
the sequence $( \phi_ {(0,\theta_2^i - \theta_1^i)})_i$ converges uniformly to $ \xi_{(0,\theta_2 - \theta_1)}$ (see \cite{Tchuia03}). 
Thus, we derive from Fathi's duality Theorem together  with the Lefschetz' property of $M$ that we have, 
 $\lim_{i\rightarrow\infty}\int_0^1[\theta_2^i(t) - \theta_1^i(t)]dt,$ is arbitrarily small. 
Therefore,
$$|[\int_0^1\theta_1^tdt - \int_0^1\theta_2^t dt]|\leq \int_0^1|(\theta_1^t - \theta_2^t) - (\theta_1^i(t) - \theta_2^i(t))|dt + 
|\int_0^1[\theta_2^i(t) - \theta_1^i(t)]dt|,$$
for all $i$. In particular, when $i$ tends to infinity, we get $|[\int_0^1\theta_1^tdt - \int_0^1\theta_2^t dt]|\leq 0$, i.e. 
 $$[\int_0^1\theta_1^tdt] =
 [\int_0^1\theta_2^t dt]. $$ This completes the proof. $\Box$

\subsubsection{Fundamental group of $SSympeo(M,\omega)$}
As far as i know, it is proved nowhere whether any continuous path in $SSympeo(M,\omega)$ is a strong symplectic isotopy or not. So, 
we cannot identify the set of all homotopic classes in  $SSympeo(M,\omega)$ generated by the loops at the the identity to be the fundamental group
 of $SSympeo(M,\omega)$. In regard of this fact, in this work, we  write $\pi_1(SSympeo(M,\omega))$ to represent 
the set of all 
the equivalent classes with respect to the equivalent relation $\sim$, whose all the representatives are loops at the identity. Therefore, 
we define a topological flux group as follows. Set $$\Gamma_0 = \widetilde{Cal_0 }(\pi_1(SSympeo(M,\omega))).$$
We have the following theorem.

\begin{theorem}\label{T7}
 Let $(M,\omega)$ be any closed symplectic manifold. Then the group $\Gamma_0$ coincides with the usual flux group $ \Gamma.$ \\

\end{theorem}

{\it Proof of Theorem \ref{T7}}. The inclusion $\Gamma \subset \Gamma_0$ is obvious.
To achieve the proof, it remains to prove that $\Gamma_0\subset\Gamma.$ For this purpose, let $[\theta]\in\Gamma_0$, by definition there exists 
$\{\xi\}\in\pi_1(SSympeo(M,\omega))$ such that $ [\theta] = \widetilde{Cal_0 }(\{\xi\} ).$ 
On the other hand, there exists a sequence $(\Phi_i)$ of symplectic isotopies such that $\xi = \lim_{C^0 + L^\infty}((\Phi_i)).$ Since $\xi$
is a loop at the identity, 
by Theorem \ref{Defor1} the sequence $(\Phi_i)$ can be deformed into a 
sequence of symplectic isotopies $(\Psi_i)$ such that $\Psi_i(1) = id$ for all $i$, and 
$\xi = \lim_{C^0 + L^\infty}((\Psi_i)).$
 More precisely, it follows from Theorem \ref{Defor1} that for each $i$, the isotopy $\Psi_i$ is of the form $ \zeta_i\circ \beta_i$ where 
$(\zeta_i)$ is some subsequence of $(\Phi_i)$, and $$\lim_{i\rightarrow\infty} Flux_\omega(\beta_i) = 0 .$$
Now, set $[\theta_i] = \widetilde{Cal_0} (\{\Psi_i\})$ for all $i$, and 
derive that for each $i$,  the element $[\theta_i]$ belongs to $ \Gamma$. 
In this way, we define a sequence $([\theta_i])_{i}$ of elements in $\Gamma$. 
The sequence $(\Psi_i\circ \xi^{-1})$ of strong symplectic isotopies converges in  $GSSympeo(M, \omega, g)$ to the identity. Then, 
we derive from the continuity of the map $ \widetilde{Cal_0}$ that $$\lim_i |\widetilde{Cal_0}(\Psi_i\circ \xi^{-1})| = 0,$$ i.e. 
$$\lim_i |[\theta_i] - [\theta]| \leq \lim_i |[\theta_i] - [\theta] - \widetilde{Cal_0}(\beta_i)| + 
\lim_i | \widetilde{Cal_0}(\beta_i)|\rightarrow0, i\rightarrow\infty,$$
since $\widetilde{Cal_0}(\Psi_i\circ \xi^{-1}) = [\theta_i] - [\theta] - \widetilde{Cal_0}(\beta_i),$ and 
$\widetilde{Cal_0}(\beta_i) = Flux_\omega(\beta_i).$
But, according to Ono \cite{Ono}, the flux group $\Gamma$ is discrete, i.e.  
$[\theta]\in \Gamma.$ Finally one concludes that the equality $\Gamma = \Gamma_0$ holds true. This completes the proof. $\Box$\\

\begin{remark} 
Theorem \ref{T7} can be viewed as a rigidity result since it suggests that the flux group $\Gamma$ stays 
invariant under the perturbation of the fundamental group $\pi_1(Symp_0(M,\omega))$ by the $(C^0 + L^\infty)-$topology at least when 
$M$ is Lefschetz.\\
\end{remark}

The homomorphism 
$ \widetilde{Cal_0}$ induces a surjective group homomorphism $Cal_0$ 
 from $SSympeo(M,\omega)$ onto $ {H}^1(M, \mathbb{R})/\Gamma$ such that the following diagram commutes 
$$
\begin{array}{lllccc}
          \widetilde{\mathcal{P}SSympeo(M, \omega)}   &\xrightarrow{\widetilde{Cal_0}} & {H}^1(M, \mathbb{R}) \\
                &                                                   &\\
    \tilde{ev} \downarrow &                                                 
 &\downarrow \pi_2\\ 
                 &                                                  &\\
    SSympeo(M,\omega) &      \xrightarrow{Cal_0}                        & {H}^1(M, \mathbb{R})/\Gamma.

\end{array}
\hspace{1cm}(II)
$$\[\]
where $\tilde{ev}$ and $\pi_2$ are projection mappings. Denote by $\widetilde{Hameo(M,\omega)}$ the subset of $ \widetilde{\mathcal{P}SSympeo(M, \omega)}$ 
which is mapped  onto the group $Hameo(M,\omega)$, under $\widetilde{ev} $. Let  $\widetilde{Ham(M,\omega)}$ 
be the universal covering of the space $Ham(M,\omega)$.
 
  \begin{theorem}\label{T8}Let $(M,\omega)$ be a closed symplectic manifold. Then, 
 $$\{\xi\}\in \widetilde{Hameo(M,\omega)}\Leftrightarrow\widetilde{Cal_0}(\{\xi\}) = 0.$$

\end{theorem}

{\it Proof of Theorem \ref{T8}}. 
Assume that $\{\xi\}\in \widetilde{\mathcal{P}SSympeo(M, \omega)}  $ such that $\{\xi\}\in \widetilde{Hameo(M,\omega)}$.  
Since $h = \xi(1)\in Hameo(M,\omega)$, then $\{\xi\}$ admits a representative 
$\lambda$ which is a continuous Hamiltonian flow. Thus,  $$\widetilde{Cal_0}(\{\xi\}) = \widetilde{Cal_0}(\{\lambda\}) = 0.$$
 For the converse, 
assume that $\{\xi\}\in \widetilde{\mathcal{P}SSympeo(M, \omega)}  $ such that $$\widetilde{Cal_0}(\{\xi\}) = 0.$$ By Fathi's Poincar\'e duality 
theorem (equation (\ref{FmF})), the Fathi mass flow of $\xi$ vanishes. Then, we derive from Theorem \ref{Defor2} that $\xi$ is homotopic to a continuous Hamiltonian 
flow $\mu$ relatively to fixed extremities, and the statement follows.$\Box$
  
 \begin{theorem}\label{T9}

 Let $\{\phi_t\}\in \widetilde{\mathcal{P}SSympeo(M, \omega)} $. Then $\phi_1$ lies in $ Hameo(M,\omega)$ if and only if  $\widetilde{Cal_0}(\{\phi_t\})$ 
lies in $ \Gamma.$ \\

\end{theorem}
{\it Proof of Theorem \ref{T9}}.  Assume that $\{\phi_t\}\in  \widetilde{\mathcal{P}SSympeo(M, \omega)} $ such that $\phi = \phi_1 \in Hameo(M,\omega)$. 
One can lift any continuous Hamiltonian flow $ \psi_t $ from the identity to $\psi_1 = \phi$ onto an 
element $\{\psi_t\}\in \widetilde{Hameo(M,\omega)}$. The strong symplectic isotopy $t\mapsto \psi_t^{-1}\circ\phi_t$ is obviously a loop at the 
identity, i.e. $ \{\psi_t^{-1}\circ\phi_t\}\in\pi_1(SSympeo(M,\omega)) $. Compute, $$\widetilde{Cal_0}(\{\psi_t^{-1}\circ\phi_t\}) = 
\widetilde{Cal_0}(\{\phi_t\}) + \widetilde{Cal_0}(\{\psi_t^{-1}\}).$$ But, we also have $\widetilde{Cal_0}(\{\psi_t^{-1}\}) =0$, 
i.e. $$\widetilde{Cal_0}(\{\phi_t\}) = 
\widetilde{Cal_0}(\{\psi_t^{-1}\circ\phi_t\})\in \Gamma.$$ For the converse, let   $\{\phi_t\}\in \widetilde{\mathcal{P}SSympeo(M, \omega)} $ such that 
$\widetilde{Cal_0}(\{\phi_t\})$ lies in $ \Gamma.$ 
Since the map $\widetilde{Cal_0}$ is surjective, there exists an element $ \{\X_t\} \in  \pi_1(SSympeo(M,\omega))$ such that 
$$\widetilde{Cal_0}(\{\phi_t\}) = \widetilde{Cal_0}(\{\X_t\}),$$ i.e. $\widetilde{Cal_0}(\{\phi_t\circ\X_t^{-1}\}) = 0.$
 Therefore, we derive from  Theorem  \ref{T8} that $\{\phi_t\circ\X_t^{-1}\}$ admits a representative which is a continuous Hamiltonian flow, 
i.e.  $ \phi_1 = \phi_1\circ\X_1^{-1} \in Hameo(M,\omega).$ 
This completes the proof. $\Box$

\subsection{On the kernel of $ Cal_0 $}

Note that a strong symplectic homeomorphism $h$ belongs to the kernel of the mapping $ Cal_0 $
 if one can find a strong symplectic isotopy $\xi_{(U,\mathcal{H})}$ 
with $\xi_{(U,\mathcal{H})}^1 = h$, and $$\widetilde{Cal_0}(\xi_{(U,\mathcal{H})}) = [\int_0^1\mathcal{H}_sds] = 0.$$ So, if 
$(\phi_{(U^i,\mathcal{H}^i)})$ is any sequence of symplectic isotopies such that 
$$\xi_{(U,\mathcal{H})} = \lim_{C^0 + L^\infty} (\phi_{(U^i,\mathcal{H}^i)}),$$ 
then one deduces from the nullity of the integral $\int_0^1\mathcal{H}_sds$ 
that for $i$ sufficiently large the flux of the isotopy $ \phi_{(U^i,\mathcal{H}^i)}$ is sufficiently small. 
Therefore we derive from Banyaga \cite{Ban78, Ban97} that the time one map 
$ \phi_{(U^i,\mathcal{H}^i)}^1$ is Hamiltonian for $i$ sufficiently large.
\begin{theorem}\label{ker!} (Topological Banyaga's theorem III) 
The  kernel of the mapping $ Cal_0 $ is path connected.\\

\end{theorem}
{\it Proof of Theorem \ref{ker!}}. Let $h$ be a strong symplectic homeomorphism that belongs 
to $ \ker Cal_0$. That is, there exists a strong symplectic isotopy $\xi_{(U,\mathcal{H})}$ 
such that $\xi_{(U,\mathcal{H})}^1 = h,$ and $\int_0^1\mathcal{H}_sds = 0.$ 
\begin{itemize}
\item Step (1). Let $\rho_t\circ \psi_t$ be the Hodge decomposition of $\xi_{(U,\mathcal{H})}$ (see \cite{Tchuia03}), 
and let $(\phi_{(U^i,\mathcal{H}^i)})$ be any sequence of symplectic isotopies such that 
$$\xi_{(U,\mathcal{H})} = \lim_{C^0 + L^\infty} (\phi_{(U^i,\mathcal{H}^i)}).$$ Then, one deduces from the 
nullity of the integral $\int_0^1\mathcal{H}_sds$ that for $i$ sufficiently large the flux 
of $ \phi_{(0,\mathcal{H}^i)}$ is sufficiently small, and therefore it follows from Banyaga \cite{Ban78, Ban97} 
that the time one map $ \phi_{(0,\mathcal{H}^i)}^1$ is Hamiltonian for $i$ sufficiently large. 
\item Step (2). Following Banyaga \cite{Ban78, Ban97} we deduce the following fact. 
Set $$\theta_i^t = \phi_{(0,\mathcal{H}^i)}^t,$$ for each $i$, and for all $t$. 
Observe that for each fixed $t$, the isotopy\\ $s\mapsto h_{(s,t)}^i = \theta_i^{st}$ connects $\theta_i^t $ 
to the identity. Now, for each $i$, consider the smooth family of harmonic vectors fields defined as follows 
: $$X_{(s,t)}^i = \frac{d }{d t} h_{(s,t)}^i\circ(h_{(s,t)}^i)^{-1},$$ and 
set $$ \alpha_t^i =  \int_0^1 i_{X_{(s,t)}^i}\omega ds =  \int_0^1\mathcal{H}^i_{st}ds.$$ It follows from the above that 
the vector field $ Y_t^i $ defined as follows $$ i_{Y_{t}^i}\omega = \alpha_t^i - t\alpha_1^i$$ is 
harmonic for each $ i $, and therefore, the following $1-$form $$ \int_0^1i_{(X_{(s,t)}^i - Y_{t}^i)}\omega ds=   t\int_0^1\mathcal{H}^i_{u}du$$ 
is harmonic. For each $ i $, let $ G_{(s,t)}^i $ be the 2-parameter family of symplectic diffeomorphisms 
defined by integrating (in $ s $) the family of harmonic vector fields $$ Z_{(s,t)}^i = X_{(s,t)}^i - Y_{t}^i.$$
 Then, Lemma \ref{LDefor2} suggests that for each fixed $t$, the 
sequence of family of symplectic vector fields $(Z_{(s,t)}^i)_s$ converges in $\|. \|^\infty$, and the sequence 
 $(G_{(s,t)}^i)_s$  of its generating paths converges uniformly. The latter limit denoted by $ (G_{(s,t)})_s $  
is then obviously a strong symplectic isotopy, and the map $ (s,t)\mapsto G_{(s,t)}$ is continuous.
 \item Step (3). Since by assumption, the quantity  $ \int_0^1 i_{(X_{(s,t)}^i - Y_{t}^i)}\omega ds$ is small for all $ i $ sufficiently large,
 it follows from Banyaga \cite{Ban78, Ban97} that the isotopy $t\mapsto G_{(1,t)}^i$ is 
Hamiltonian for all  sufficiently large $ i $. That is, the limit $ t\mapsto G_{(1,t)} $ is a 
continuous Hamiltonian flow, i.e. a continuous path in $ \ker Cal_0$. 
 For instance, we compute $ G_{(1,1)} = \lim_{C^0}(\theta_i^1) = \rho_1,$ and deduce that 
the mapping  $ t\mapsto G_{(1,t)}\circ\psi_t $ is a  continuous path in $ \ker Cal_0$ 
that connects $ G_{(1,1)}\circ\psi_1 = h $ to the identity. This completes the proof. $ \Box$\\
 \end{itemize}
 \begin{theorem}\label{ker!2}  Any strong symplectic isotopy in $ \ker Cal_0$ is a continuous Hamiltonian flow.
\end{theorem}
{\it Proof of Theorem \ref{ker!2}}. We will adapt a proof given in Banyaga \cite{Ban78, Ban97} into our case. 
\begin{itemize}
\item Step (1). Let $ \sigma_{(U,\mathcal{H})} $ be a strong symplectic isotopy 
in $ \ker Cal_0$, and let  $(\phi_{(U^i,\mathcal{H}^i)})$ be any sequence of symplectic isotopies such that 
$$ \sigma_{(U,\mathcal{H})} = \lim_{C^0 + L^\infty} (\phi_{(U^i,\mathcal{H}^i)}).$$ Since 
$\pi_1(SSympeo(M,\omega))$ acts on $ \widetilde{\mathcal{P}SSympeo(M, \omega)}  $, we denote by
\\ $\pi_1(SSympeo(M,\omega)). \ker \widetilde{Cal_0}$ the set of all the 
orbits of the points of $\ker \widetilde{Cal_0}$. 
It is clear that the mapping $$ \tilde{ev} : \widetilde{\mathcal{P}SSympeo(M, \omega)}  \rightarrow SSympeo(M,\omega),$$ is the covering map. Then  we have 
$$ \tilde{ev}^{-1}( \ker Cal_0) = \pi_1(SSympeo(M,\omega)). \ker \widetilde{Cal_0}. $$ Let $ \widetilde\sigma $ be the 
lifting of the path $ \sigma_{(U,\mathcal{H})} $ such that $ \widetilde\sigma(0) = id_M $. By assumption,  
we have $$\widetilde{Cal_0}(\widetilde\sigma(t))\in \Gamma,$$ for all $t$. A verbatim repetition of some 
arguments from Banyaga \cite{Ban78, Ban97} which are supported by the discreteness of the flux group $\Gamma$ 
leads to $$ \widetilde\sigma(t) \in \ker \widetilde{Cal_0},$$ for all $t$. We then derive from the 
lines of the proof of Theorem \ref{ker!} that for each $t$, the homotopy class $  \widetilde\sigma(t)$ 
admits a representative $ s\mapsto G_{s,t}$  which is  a continuous Hamiltonian flow, i.e.
 a continuous path in $ \ker Cal_0$, and the map $ (s, t)\mapsto G_{s,t}$ is continuous. For all fixed $ t $, the following map
: $$ (u, s)\mapsto G_{(u(s-1) + 1),((u-1)s + u)t} $$
 induces a homotopy between $ s\mapsto\sigma_{(U,\mathcal{H})}^{(t.s)} $ and $ s\mapsto G_{s,t}$.  

\item Step (2). It follows from step (1) that $$ \widetilde{Cal_0}(\sigma(t.s)) =  \int_0^t\mathcal{H}_u du  = 0 ,$$ for all $t$. On the other hand, since 
$$ \sigma_{(U,\mathcal{H})} = \lim_{C^0 + L^\infty} (\phi_{(U^i,\mathcal{H}^i)}),$$ we derive 
 that the strong symplectic isotopy $$ \Psi_t : s\mapsto \sigma_{(U,\mathcal{H})}^{(s.t)},$$ is the $ (C^0 + L^\infty)-$limit 
of the sequence of symplectic isotopies $$ \Psi_t^i : s\mapsto \phi_{(U^i,\mathcal{H}^i)}^{(s.t)}.$$ For each $t$, the symplectic 
isotopy $ \Psi_t^i : s\mapsto \phi_{(U^i,\mathcal{H}^i)}^{(s.t)} $ is generated by 
$  (V_t^i, \mathcal{K}_t^i)$ where $ V_t^i(s) = tU_{(st)}^i,$ and $\mathcal{K}_t(s)^i = t\mathcal{H}_{(st)}^i$ for all $s$. 
Therefore, the fact that 
$$ \widetilde{Cal_0}(\sigma(t.s)) = 0,$$ for all $t$ suggests that the flux of the isotopy $\Psi_t^i$ tends to 
zero when $i$ goes to the infinity, i.e. $ \int_0^t\mathcal{H}^i_u du  \rightarrow0 , i\rightarrow \infty$ for all $t$. 
That is, 
$$\dfrac{1}{h}(\int_0^{t + h}\mathcal{H}^i_u du - \int_0^t\mathcal{H}^i_u du)\rightarrow0 , i\rightarrow \infty,$$
for all $h$ smaller and positive such that $(h + t)\in[0,1]$, i.e. 
$$\lim_{i\rightarrow\infty}(\lim_{h\rightarrow 0^+}\dfrac{1}{h}\int_t^{t + h}\mathcal{H}^i_u du) = 
\lim_{i\rightarrow \infty} (\mathcal{H}^i_t) =
0,$$
for all $t\in [0,1[,$ and since for each $i$ the function $ t\mapsto \mathcal{H}^i_t$ is continuous, we derive that 
$\lim_{i\rightarrow \infty} (\mathcal{H}^i_t) = 0,$ for all $t\in [0,1].$ 
This implies that $\sup _t |\mathcal{H}^i_t |\rightarrow 0, i\rightarrow \infty,$ i.e. 
$$ \sup _t |\mathcal{H}_t | \leq \sup _t |\mathcal{H}^i_t - \mathcal{H}_t| + \sup _t |\mathcal{H}^i_t |\rightarrow 0, i\rightarrow \infty.$$ 
Finally, we see that the harmonic part $\mathcal{H}$ of the  generator of the $ (U,\mathcal{H}) $ is trivial, 
i.e. $\sigma_{(U,\mathcal{H})} = \sigma_{(U,0)}.$
That is, the strong symplectic isotopy $ \sigma_{(U,\mathcal{H})} $ is in fact a continuous Hamiltonian flow. $\Box$
\end{itemize}
 \begin{proposition}\label{hameo 1} $
 Hameo(M,\omega)$ coincides with $ \ker Cal_0$. 
 
 \end{proposition}

 {\it Proof of Proposition  \ref{hameo 1}}. According to Oh-M\"{u}ller \cite{Oh-M07}, 
any Hamiltonian homeomorphism $h$ can be connected to the identity trough a continuous Hamiltonian flow $\xi$,
 we then derive that $ \widetilde{Cal_0}(\{\xi\}) = 0,$ which implies $ Cal_0(h) = 0.$
 The converse inclusion follows from  Theorem \ref{T9}. In fact,  
$ h\in \ker Cal_0$ implies that for any equivalent class such that $\xi(1) = h,$
 we have $\pi_2(\widetilde{Cal_0}({\xi})) = 0,$ i.e. $ \widetilde{Cal_0}({\xi})\in \Gamma$. 
This is equivalent to say that $ h\in Hameo (M,\omega)$. $\Box$\\
 
\begin{theorem}\label{EN1311}(Topological Banyaga's theorem I) Let $(M,\omega)$ be any closed symplectic manifold. Then,  
 any strong symplectic isotopy in the group of Hamiltonian homeomorphisms is a continuous Hamiltonian flow. 
\end{theorem}
{\it Proof}. The proof of Theorem \ref{EN1311} is 
 a verbatim repetition of the proof of Theorem \ref{ker!2}. $\Box$ 
\begin{theorem}\label{EN131}(Topological Banyaga's theorem II) Let $(M,\omega)$ be any closed symplectic manifold. 
 Then, the group of Hamiltonian homeomorphisms is path connected, and locally connected.
 
\end{theorem}

Note that Oh-M\"{u}ller \cite{Oh-M07} proved that the group of Hamiltonian homeomorphisms is locally path 
connected, then locally connected. But, in the present paper we give a different proof.\\

{\it Proof of Theorem \ref{EN131}}.
 The result of Theorem \ref{ker!} states that $Hameo(M,\omega) = \ker Cal_0,$ 
is path connected, while the result of Theorem \ref{T7} implies that  $\Gamma = \widetilde{Cal_0}(\pi_1(SSympeo(M,\omega))$ is discrete. 
Thus, $Hameo(M,\omega) = \ker Cal_0$ is locally connected .$\Box$ 
\begin{theorem}\label{EN13}(Weak topological Weinstein's theorem) Let $(M,\omega)$ be a closed symplectic manifold. 
 Then, the group of strong symplectic homeomorphisms is locally path connected.
 \end{theorem}

 {\it Proof}.
Since $ SSympeo(M,\omega) $ is a topological group, 
the latter is a homogeneous space, i.e. for all $p, q \in SSympeo(M,\omega)$, there exists 
a homeomorphism $ \Phi $ that maps $p$ onto $q$, i.e. $ \Phi(p) = q .$ 
So, it suffices just to check a local property at a point (namely, at the identity) to prove it on the entire group. 
But, the result of Theorem \ref{ker!2} states that any strong 
symplectic isotopy in $Hameo(M,\omega)$ is a continuous Hamiltonian flow, while a result dues 
to Oh-M\"{u}ller \cite{Oh-M07} states that $ Hameo(M,\omega) $ is locally path connected. 
Therefore, one derives from the above statements that for every open neighborhood  
$V\subseteq SSympeo(M,\omega) $ of the constant map identity, there exists a path connected 
open set $U\subset Hameo(M,\omega)$ with $id \in U \subset V$. Now, let $p\in SSympeo(M,\omega) $, 
and let $V_p\subseteq SSympeo(M,\omega) $ be an arbitrary neighborhood of $p$. There exists a 
homeomorphism $ \Phi $ that sends the identity onto $p$, i.e. $ \Phi(id) = p .$ So, $ \Phi$ displaces 
some neighborhood $V_0$ of the constant map identity ( with $V_0\subset U$ ) onto a neighborhood $U_0$ of $ p $. 
But, according to Oh-M\"{u}ller \cite{Oh-M07}, $V_0$ contains a path connected neighborhood of the identity, 
and the latter is displaced by $\Phi$ onto a path connected neighborhood of $p$ denoted by $O_p$ with $O_p\subset U_0$. 
That is, $  O_p\cap V_p$ is a path connected neighborhood of $p$ which is contained in $ V_p$. This completes the proof. $\Box$ 

\section{$C^0-$Banyaga's Hofer-like norm}
Actually, we focus on the construction of a $C^0-$Hofer-like norm for
 strong symplectic homeomorphisms. In fact the following results points out the existence and some studies of 
the $C^0-$analogue of the well-known Banyaga's Hofer-like norm found in \cite{Ban08a}. 
We recall that the symplectic topology on $SSympeo(M, \omega)$ is defined to be the strongest topology on it which makes  
the map $$ev:\mathcal{P}SSympeo(M, \omega)\rightarrow SSympeo(M, \omega),$$ $$\gamma_{(U,\mathcal{H})}\mapsto \gamma_{(U,\mathcal{H})}^1,$$ 
becomes a surjective homomorphism in 
the category of topological groups with respect to the $ (C^0 + L^\infty)-$topology on the space $\mathcal{P}SSympeo(M, \omega)$. 
\begin{definition}
 For any strong symplectic isotopy $\gamma_{(U,\mathcal{H})}$, we define its length by  
$$l(\gamma_{(U,\mathcal{H})}) = \dfrac{l_\infty(\gamma_{(U,\mathcal{H})}) + l_\infty(\gamma_{(U,\mathcal{H})}^{-1})}{2},$$
where $l_\infty(\gamma_{(U,\mathcal{H})}) = \max_t(osc(U_t) + |\mathcal{H}_t|).$

\end{definition}
\begin{definition}
 For any $h\in SSympeo(M,\omega)$, we define the Hofer-like norm of $h$ by 
\begin{equation}\label{SHL}
 \|h\|_{SHL} = : \inf\{l(\gamma_{(U,\mathcal{H})})\},
\end{equation}
 where the infimum is taken over all strong symplectic isotopies $\gamma_{(U,\mathcal{H})}$ with 
$\gamma_{(U,\mathcal{H})}^1 = h$.
\end{definition}
The above definitions make sense because of the uniqueness results (Theorem \ref{UG} and Theorem \ref{URBT}). 
This agrees with the usual definition of Hofer-like length of symplectic isotopies \cite{Ban08a}.\\

\begin{theorem}\label{SHLT1} (Topological Banyaga's norm ) Let $(M,\omega)$ be any closed symplectic manifold. 
 Then, the mapping  $$\|.\|_{SHL} : \mathfrak{SS}ympeo(M,\omega)\rightarrow \mathbb{R},$$ 
is continuous.
\end{theorem}
{\it Proof.} 
We will process step by step.
\begin{itemize}
 \item Step (1). Positivity is obvious. 
In the way that the symplectic topologies has been defined on the spaces $SSympeo(M,\omega)$ and $GSSympeo(M, \omega, g)$ (see \cite{BanTchu}), 
the map $$\mathfrak{SS}ympeo(M,\omega)\rightarrow \mathbb{R},$$ $$h\mapsto\|h\|_{SHL},$$ is continuous if and only if 
$$ \mathcal{P}SSympeo(M, \omega)\xrightarrow{ev} \mathfrak{SS}ympeo(M,\omega)\xrightarrow{\|.\|_{SHL}} \mathbb{R},$$ 
is continuous with respect to the symplectic topology on the space\\ $\mathcal{P}SSympeo(M, \omega)$. 
 Then, the continuity of 
the map $$\mathfrak{SS}ympeo(M,\omega)\rightarrow \mathbb{R},$$ $$h\mapsto\|h\|_{SHL},$$ is proved as follows. Let 
$\gamma_{(U_k,\mathcal{H}_k)}$ be a sequence of elements of $\mathcal{P}SSympeo(M, \omega)$ 
that converges to $\gamma_{(U,\mathcal{H})}\in \mathcal{P}SSympeo(M, \omega)$ with respect to the $(C^0 + L^\infty)-$topology. We always have:
$$2|\|\gamma(1)\|_{SHL} - \|\gamma_k(1)\|_{SHL}|\leq 2\|(\gamma_{(U_k,\mathcal{H}_k)}^{-1}\circ \gamma_{(U,\mathcal{H})})(1)\|_{SHL}$$
$$\leq 2l[\gamma_{(U_k,\mathcal{H}_k)}^{-1}\circ \gamma_{(U,\mathcal{H})}],$$
$$\leq \max_{t\in[0,1]}\{osc(U_k^t\circ\gamma_{(U_k,\mathcal{H}_k)} -U^t\circ\gamma_{(U_k,\mathcal{H}_k)}) + |\mathcal{H}^t - \mathcal{H}_k^t|\}$$ $$
 + \max_{t\in[0,1]}osc(\Delta^0_t(\mathcal{H}_k  ,\gamma_{(U_k,\mathcal{H}_k)}) -\Delta^0_t( \mathcal{H} ,\gamma_{(U_k,\mathcal{H}_k)}))$$
$$+ \max_{t\in[0,1]}\{osc(U_k^t\circ\gamma_{(U,\mathcal{H})} -U^t\circ\gamma_{(U,\mathcal{H})}) + |\mathcal{H}^t - \mathcal{H}_k^t|\} 
$$ $$+ \max_{t\in[0,1]} osc(\Delta^0_t(\mathcal{H}, \gamma_{(U,\mathcal{H})}) -\Delta^0_t( \mathcal{H}_k ,\gamma_{(U,\mathcal{H})})).$$ 
\item Step (2). By assumption the sequence $\gamma_{(U_k,\mathcal{H}_k)}$  converges to $\gamma_{(U,\mathcal{H})}$ 
 with respect to the $(C^0 + L^\infty)-$topology, as well as the sequence
$\gamma_{(U_k,\mathcal{H}_k)}^{-1}$  converges to $\gamma_{(U,\mathcal{H})}^{-1}$ 
 with respect to the $(C^0 + L^\infty)-$topology 
since $\mathcal{P}SSympeo(M, \omega)$ is a topological group 
with respect to the symplectic topology. Therefore we derive from the above statement that 
$$\max_{t\in[0,1]}\{osc(U_k^t\circ\gamma_{(U_k,\mathcal{H}_k)} -U^t\circ\gamma_{(U_k,\mathcal{H}_k)} ) + |\mathcal{H}^t - \mathcal{H}_k^t|\} 
= \max_{t\in[0,1]}\{osc(U_k^t -U^t) + |\mathcal{H}^t - \mathcal{H}_k^t|\},$$
$$\max_{t\in[0,1]}\{osc(U_k^t\circ\gamma_{(U,\mathcal{H})} -U^t\circ\gamma_{(U,\mathcal{H})}) + |\mathcal{H}^t - \mathcal{H}_k^t|\} = 
\max_{t\in[0,1]}\{osc(U_k^t -U^t) + |\mathcal{H}^t - \mathcal{H}_k^t|\},
$$
where the right-hand side of the above estimates tends to zero when $i$ goes at infinity. 
\item Step (3). Next, observe that
$$\max_{t\in[0,1]}osc(\Delta^0_t(\mathcal{H}  ,\gamma_{(U,\mathcal{H})}) -\Delta^0_t( \mathcal{H}_k ,\gamma_{(U,\mathcal{H})})) 
= \max_{t\in[0,1]}osc(\Delta^0_t(\mathcal{H}_k - \mathcal{H},\gamma_{(U,\mathcal{H})}),$$
$$\max_{t\in[0,1]}osc(\Delta^0_t(\mathcal{H}_k, \gamma_{(U_k,\mathcal{H}_k)}) -\Delta^0_t( \mathcal{H}, \gamma_{(U_k,\mathcal{H}_k)})) 
= \max_{t\in[0,1]}osc(\Delta^0_t(\mathcal{H}_k - \mathcal{H}, \gamma_{(U_k,\mathcal{H}_k)}).$$
 On the other hand, by definition of $\mathcal{P}SSympeo(M, \omega)$,  
for any fixed integer $ k$, there exists a sequence of symplectic isotopies $(\phi_{(V^{k, j}, \mathcal{K}^{k, j})})$   such that 
$$ \gamma_{(U_k,\mathcal{H}_k)} = \lim_{C^0 + L^\infty} ((\phi_{(V^{k, j}, \mathcal{K}^{k, j})}).$$
 Similarly, 
 there exists a  sequence of symplectic isotopies $(\phi_{(U^i, \mathcal{H}^i)})$
  which does not depend on $k$  such 
that $$ \gamma_{(U,\mathcal{H})} = \lim_{C^0 + L^\infty} (\phi_{(U^i, \mathcal{H}^i)}).$$
We compute, 
 $$\Delta^0(\mathcal{H}_k - \mathcal{H},\gamma_{(U,\mathcal{H})}) = \lim_{L^\infty}\{\int_0^t(\mathcal{K}^{k, j}_t - 
\mathcal{H}^i_t)(\dot\phi_{(U^i, \mathcal{H}^i)}^s )\circ\phi_{(U^i, \mathcal{H}^i)}^s ds\},$$
 $$\Delta^0(\mathcal{H} - \mathcal{H}_k,\gamma_{(U_k,\mathcal{H}_k)}) = \lim_{L^\infty}\{\int_0^t( 
\mathcal{H}^i_t -\mathcal{K}^{k, j}_t)(\dot\phi_{(V^{k, j}, \mathcal{K}^{k, j})}^s )\circ\phi_{(V^{k, j}, \mathcal{K}^{k, j})}^s ds\}.$$
(see \cite{BanTchui02} for more convenience).
\item Step (4). Since both sequences $(\phi_{(U^i, \mathcal{H}^i)})_i$ and $(\phi_{(V^{k, j}, \mathcal{K}^{k, j})})_j$ are Cauchy in $\bar d$, 
and we have $$\max_{t\in[0,1]}| \mathcal{H}^j_t -\mathcal{K}^{k, j}_t|\rightarrow0,$$ when $j\textgreater k$ and $k\rightarrow \infty$; 
we derive from \cite{BanTchu} (Lemma 3.4, \cite{BanTchu}) that, 
$$\max_{t\in[0,1]}osc(\Delta^0_t(\mathcal{H}_k - \mathcal{H},\gamma_{(U_k,\mathcal{H}_k)} )\rightarrow0, k\rightarrow\infty,$$ 
$$\max_{t\in[0,1]}osc(\Delta^0_t(\mathcal{H} - \mathcal{H}_k,    \gamma_{(U,\mathcal{H})})\rightarrow0, k\rightarrow\infty.$$
\item Step (5). 
 Finally, the results of steps (2) and (4) show that 
$$l[\gamma_{(U_k,\mathcal{H}_k)}^{-1}\circ \gamma_{(U,\mathcal{H})}] \rightarrow0, k\rightarrow\infty.$$ This completes the proof. $\Box$\\

\end{itemize}
Here, we point out some properties of the mapping $\|.\|_{SHL}.$
\begin{theorem}\label{SHLT2} Let $(M,\omega)$ be any closed symplectic manifold. 
Given two strong symplectic homeomorphisms $ h,$ and $f$ we have : 
\begin{enumerate}
\item $\|h\|_{SHL} = 0$ if $h = id$,
\item $\|h\|_{SHL} = \|h^{-1}\|_{SHL}$,
\item $\|h\circ f\|_{SHL}\leq \|h\|_{SHL} + \|f\|_{SHL}$.\\            
\end{enumerate}

\end{theorem}
{\it Proof.} The item (2) follows from the definition of the map $\|,\|_{SHL}$. For (1),
 we adapt the proof of the nondegeneracy of the Hofer norm for 
Hamiltonian homeomorphism given by Oh \cite{Oh}. We will proceed step by step.
\begin{itemize}
 \item Step (a). Suppose that $h\neq id$. 
Then $h$ displaces a small nonempty compact ball $B$ of positive symplectic displacement energy 
$e_S(B)\textgreater 0$. For such a ball $B$, we set $\delta = e_S(B)\textgreater 0.$
From characterization of the infimum, one can find a strong symplectic isotopy $\gamma_{(U,\mathcal{H})}$  with $\gamma_{(U,\mathcal{H})}^1 = h$, 
such that $\|h\|_{SHL}\textgreater l(\gamma_{(U,\mathcal{H})}) - \dfrac{\delta}{4}.$
\item Step (b). On the other hand, it follows from the definition of the topological group  
$GSSympeo(M, \omega, g)$ that   
there exists a sequence $(\phi_{(F_i,\lambda_i)})$ that converges to $\gamma_{(U,\mathcal{H})}$ with respect to the 
$(C^0 + L^\infty)-$topology. 
So, we can find an integer $i_0$ with $(\phi_{(F_{i_0},\lambda_{i_0})})$ 
 is sufficiently close to $\gamma_{(U,\mathcal{H})}$ [resp. $(\phi_{(F_{i_0},\lambda_{i_0})}^{-1})$ 
 sufficiently close to $\gamma_{(U,\mathcal{H})}^{-1}$] with respect to the 
$(C^0 + L^\infty)-$topology so that 
$$l(\gamma_{(U,\mathcal{H})}) \textgreater l(\phi_{(F_{i_0},\lambda_{i_0})}) - \dfrac{\delta}{4},$$ where $\phi = 
\phi_{(F_{i_0},\lambda_{i_0})}^1$ displaces $B$. 
\item Step (c). It follows from the definition of Banyaga's Hofer-like norm $\|,\|_{HL} $ that 
$$l(\phi_{(F_{i_0},\lambda_{i_0})}) - \dfrac{2\delta}{4}
\geq \|\phi\|_{HL} - \dfrac{\delta}{2},$$ and by definition of the symplectic displacement energy $e_S$, we have 
$$ \|\phi\|_{HL} - \dfrac{\delta}{2}\geq \delta - \dfrac{\delta}{2} = \dfrac{\delta}{2}\textgreater 0.$$ 
\item Step (d). The statements of steps (a), (b) and (c) imply that 
$$\|h\|_{SHL}\textgreater l(\gamma_{(U,\mathcal{H})}) - \dfrac{\delta}{4} \textgreater 
l(\phi_{(F_{i_0},\lambda_{i_0})}) - \dfrac{2\delta}{4}
\geq \|\phi\|_{HL} - \dfrac{\delta}{2}\textgreater 0. $$ 
 The item (3) follows from the continuity of the map $\|,\|_{SHL}$ closely the proof triangle inequality of Banyaga's Hofer-like norm.
 This completes the proof. $\Box$\\ 
\end{itemize}

Theorem \ref{SHLT2} and Theorem \ref{SHLT1} are similar to some results that was proved by Oh \cite{Oh} 
for Hamiltonian homeomorphisms. In particular, if the manifold $M$ is simply connected, 
then the norm $\|.\|_{SHL}$ reduces to the norm $\|.\|_{Oh}$ constructed by Oh on the space of 
all Hamiltonian homeomorphisms \cite{Oh-M07, Oh}. Otherwise, when $H^1(M,\mathbb{R}) \neq 0$,   
it follows from the decomposition theorem of strong symplectic homeomorphisms that
 the group of Hamiltonian homeomorphisms is strictly contained in the group of 
strong 
symplectic homeomorphisms (see \cite{Tchuia03}). So, it is important to know whether the norm $\|.\|_{SHL}$ is an extension of Oh's norm or not. 
This leads to the following result.
\begin{theorem}\label{EN}(Topological Banyaga's conjecture)
 The norm $\|.\|_{SHL}$ restricted to the group of Hamiltonian homeomorphisms is equivalent to Oh's norm.
\end{theorem}
This result is motivated in party by a conjecture from Banyaga \cite{Ban08a} that was proved by Buss-Leclercq \cite{BusLec11}. 
The latter conjecture 
stated that the restriction of Banyaga's Hofer-like norm to the group of Hamiltonian 
 diffeomorphisms is equivalent to Hofer's norm. \\

{\it Proof of Theorem \ref{EN}.} By construction, we always have $$\|.\|_{SHL}\leq \|.\|_{Oh}.$$ 
To complete the proof, we need to show that there exists a positive finite constant 
$A_0$ such that $$\|.\|_{Oh}\leq A_0.\|.\|_{SHL},$$ or equivalently, via the sequential criterion, it suffices to prove that
any sequence of Hamiltonian homeomorphisms converging to the constant map identity
for  $\|.\|_{SHL}$, converges to the constant map identity for Oh's norm $\|.\|_{Oh}$. 
The proof that we give here heavily relies the ideas that Oh \cite{Oh} and Banyaga \cite{Ban08a} 
used in the proof of the nondegeneracy of their norms. Let $\psi^i$ be a sequence of Hamiltonian homeomorphisms that converges to the identity 
with respect to the norm $\|.\|_{SHL}$. For each $i$, and any 
$\epsilon \textgreater0$ there exists a strong symplectic isotopy
 $\psi_{(U^{i,\epsilon},\mathcal{H}^{i,\epsilon})}$ such that  
 $ \psi_{(U^{i,\epsilon},\mathcal{H}^{i,\epsilon})}^1 = \psi^i,$ and 
$l(\psi_{(U^{i,\epsilon},\mathcal{H}^{i,\epsilon})})\textless \|\psi^i\|_{SHL} + \epsilon.$ 
On the other hand, for a fixed $i$, there exists a sequence of symplectic isotopies 
$(\phi_{(V_{i, j}, \mathcal{K}_{i, j})})$  such that 
$$ \psi_{(U^{i,\epsilon},\mathcal{H}^{i,\epsilon})} = \lim_{C^0 + L^\infty} (\phi_{(V_{i, j}, \mathcal{K}_{i, j})}).$$
  In particular, one can find 
an integer $j_0$ with $\phi_{(V_{i,j_0}, \mathcal{K}_{i,j_0})}$ 
sufficiently close to $\psi_{(U^{i,\epsilon},\mathcal{H}^{i,\epsilon})}$ with respect to the 
$(C^0 + L^\infty)-$topology so that 
$$l(\psi_{(U^{i,\epsilon},\mathcal{H}^{i,\epsilon})})\textgreater 
l(\phi_{(V_{i,j_0}, \mathcal{K}_{i,j_0})}) - \dfrac{\epsilon}{4}.$$ Observe that 
when i tends to infinity, Banyaga's length of the symplectic isotopy $\phi_{(V_{i,j_0}, \mathcal{K}_{i,j_0})}$ 
can be considered as being sufficiently small. Thus, 
it comes from Banyaga \cite{Ban08a} that the time-one map of such an isotopy is a Hamiltonian diffeomorphism. 
So, we can assume (without breaking the generality) that 
$\phi_{(V_{i,j_0}, \mathcal{K}_{i,j_0})}^1$ is Hamiltonian for $i$ sufficiently large. Now, we derive from the above statements that 
$$l(\psi_{(U^{i,\epsilon},\mathcal{H}^{i,\epsilon})})\textgreater
 \|\phi_{(V_{i,j_0}, \mathcal{K}_{i,j_0})}^1\|_{HL} - \dfrac{\epsilon}{4}.$$
  For $i$ sufficiently large, since the diffeomorphism $\phi_{(V_{i,j_0}, \mathcal{K}_{i,j_0})}^1$ is Hamiltonian, it 
follows from a result from Buss-Leclercq \cite{BusLec11} that there exists a positive finite constant $D$ which does not depend on $i$ such that 
$$\dfrac{1}{D}\|\phi_{(V_{i,j_0}, \mathcal{K}_{i,j_0})}^1\|_{H}\leq \|\phi_{(V_{i,j_0}, \mathcal{K}_{i,j_0})}^1\|_{SHL},$$ where 
$\|.\|_{H}$ represents the Hofer norm of Hamiltonian diffeomorphisms. The above statements imply that for $i$ sufficiently large,
$$\|\psi^i\|_{SHL} + \epsilon \textgreater \dfrac{1}{D}\|\phi_{(V_{i,j_0}, \mathcal{K}_{i,j_0})}^1\|_{H}- \dfrac{\epsilon}{4}.$$ At this level, 
we use the fact that Oh's norm restricted to the group of Hamiltonian 
diffeomorphisms is bounded from above by Hofer's norm to get 
$$ \|\psi^i\|_{SHL} + \epsilon \textgreater  \dfrac{1}{D}\|\phi_{(V_{i,j_0}, 
\mathcal{K}_{i,j_0})}^1\|_{Oh}- \dfrac{\epsilon}{4},$$ for $i$ sufficiently large. Passing to the limit in the latter estimate, one obtains 
$$ \lim_{i\rightarrow\infty}\|\psi^i\|_{SHL} + \dfrac{5\epsilon}{4}\geq 
\dfrac{1}{D}\lim_{j_0\geq i\hspace{0.cm}, i\rightarrow\infty}\|\phi_{(V_{i,j_0}, \mathcal{K}_{i,j_0})}^1\|_{Oh}= 
\dfrac{1}{D}\lim_{i\rightarrow\infty}\|\psi^i\|_{Oh},$$ for all $\epsilon$. 
Finally, we have proved that for all positive real number $\delta$ (replacing $\epsilon$ by $\dfrac{4\delta}{5D}$), we have
$$ \delta\geq \lim_{i\rightarrow\infty}\|\psi^i\|_{Oh}.$$ This completes the proof.$\square$ 
\section{Topological symplectic displacement energy}
\begin{definition}
 The strong symplectic displacement energy $ e_S^0(B)$ of a non empty compact set $B\subset M$ is :
$$e^0_S(B) = \inf\{\|h\|_{SHL} |  h\in SSympeo(M,\omega), h(B)\cap B = \emptyset\}.$$
\end{definition}
\begin{lemma}\label{SSE}
 For any non empty compact set $B\subset M$, $ e_S^0(B)$ is a strict positive number.
\end{lemma}
{\it Proof of Lemma \ref{SSE}.} 
Let $\epsilon \textgreater 0$. By definition of $ e_S^0(B)$ there exists a strong symplectic isotopy 
$\psi_{(F^{\epsilon},\lambda^{\epsilon})}$ such that 
 $\psi_{(F^{\epsilon},\lambda^{\epsilon})}^1 = h,$ and $e^0_S(B)\textgreater\|h\|_{SHL} - \epsilon.$ On the other hand, 
there exists a sequence 
of symplectic isotopies $(\phi_{(F_i,\lambda_i)})$ that converges 
to $\psi_{(F^{\epsilon},\lambda^{\epsilon})}$ whit respect to the $(C^0 + L^\infty)-$topology so that  $\phi_{(F_i,\lambda_i)}^1$ 
 displaces $B$ for $i$ sufficiently large. It follows from Theorem \ref{BDS} that 
$\|\phi_{(F_i,\lambda_i)}^1\|_{HL}\geq e_S(B)\textgreater 0$ for $i$ sufficiently large. 
Since $\epsilon \textgreater 0$ is arbitrary, the continuity of the map $\|.\|_{SHL}$  imposes that 
$$\|h\|_{SHL} - \epsilon\geq \|\phi_{(F_i,\lambda_i)}^1\|_{SHL} = \|\phi_{(F_i,\lambda_i)}^1\|_{HL},$$ for $i$ sufficiently large. Therefore,   
 $$e^0_S(B)\textgreater\|h\|_{SHL} - \epsilon\geq e_S(B)\textgreater 0. \Box$$
 
 \subsection{$ L^{(1, \infty)}-$norm for strong symplectic isotopies }
 For any strong symplectic isotopy $\gamma_{(U,\mathcal{H})}$, we define its interpolation length by  
$$l^{(1, \infty)}(\gamma_{(U,\mathcal{H})}) = \dfrac{l_0(\gamma_{(U,\mathcal{H})}) + l_0(\gamma_{(U,\mathcal{H})}^{-1})}{2},$$
with $l_0(\gamma_{(U,\mathcal{H})}) = \int_0^1(osc(U_t) + |\mathcal{H}_t|)dt.$ Therefore, for any $h\in SSympeo(M,\omega)$, 
we define the $L^{(1, \infty)}$ Hofer-like norm of $h$ by 
\begin{equation}\label{SHL}
 \|h\|_{SHL}^{(1, \infty)} = : \inf\{l^{(1, \infty)}(\gamma_{(U,\mathcal{H})})\},
\end{equation}
 where the infimum is taken over all strong symplectic isotopies $\gamma_{(U,\mathcal{H})}$ with 
$\gamma_{(U,\mathcal{H})}^1 = h$.

\subsubsection*{Conjecture A}

 Let $(M,\omega)$ be a closed connected symplectic manifold. For any $h\in SSympeo(M,\omega)$, 
we have,
 $$\|h\|_{SHL} =  \|h\|_{SHL}^{(1, \infty)},$$ holds true.\\

Conjecture A is supported by the uniqueness result of Hofer-like geometry from \cite{TD}.

\subsubsection*{Conjecture B} 
 Let $(M,\omega)$ be a closed connected symplectic manifold. 
Let $h\in SSympeo(M,\omega)$. The norm $\phi \mapsto \|h\circ\phi \circ h^{-1}\|_{SHL}$ is equivalent to the norm
$\phi \mapsto \|\phi\|_{SHL}$.\\

Conjecture B is supported by a result from \cite{BDS} (Theorem 7).

\section{Examples} 
 
In the following section, we provide some examples of the above constructions.
\subsubsection*{Example 6.1} 

 We know that any smooth symplectic isotopy is a trivial strong symplectic isotopy.
\subsubsection*{Example 6.2} 

 According to Theorem \ref{UGR}, any strong symplectic isotopy which is a $1-$parameter group  
admits a time independent generator. A harmonic 1-parameter group is an isotopy $\beta = \{\beta_t\}$ generated by the vector field 
$X$ defined by $\iota(X)\omega = \mathcal{K},$ where $\mathcal{K}$
is a harmonic 1-form. For instance, let $\phi$ be a Hamiltonian homeomorphism. Then, by definition of $\phi$, 
there exists a sequence of Hamiltonian isotopies $ \Phi_j $ which is Cauchy in $(C^0 + L^\infty)$ such that $\phi_j := \Phi_j(1) \rightarrow \phi$ 
uniformly. The sequence of symplectic isotopies $\Psi_j : t \mapsto \phi_j^{-1} \circ \beta_t\circ \phi_j,$ connects the identity to 
$ \phi_j^{-1} \circ \beta_1\circ \phi_j$. The latter sequence of symplectic isotopies is generated by 
$(\int_0^1\mathcal{K}(\dot{\Phi}_j(s))\circ\Phi_j(s)ds , \mathcal{K})$.  
In fact, the Hofer's norm of the smooth function $$ x\mapsto \int_0^1\mathcal{K}(\dot{\Phi}_j(s))\circ\Phi_j(s)ds(x),$$ does not 
depend on the choice of the isotopy that connects $\phi_j := \Phi_j(1)$ to the identity, since for each $j$, we always 
have, $$\phi_j^\ast (\mathcal{K}) - \mathcal{K} =   d(\int_0^1\mathcal{K}(\dot{\alpha}_j(s))\circ\alpha_j(s)ds),$$ for all isotopy 
$\alpha_j$ that connects 
$\phi_j := \Phi_j(1)$ to the identity ( see \cite{BanTchu}, Definition 2.2 for more details). As any 
reader can see, the sequence of generators 
$(\int_0^1\mathcal{K}(\dot{\Phi}_j(s))\circ\Phi_j(s)ds , \mathcal{K})$ is Cauchy in $D^2$ if and only if the sequence of functions
$ x\mapsto (\int_0^1\mathcal{K}(\dot{\Phi}_j(s))\circ\Phi_j(s)ds)(x)$ is Cauchy in the $L^\infty$ Hofer norm. But, 
Lemma 3.9 from  \cite{TD} shows that the sequence of functions
$$ x\mapsto (\int_0^1\mathcal{K}(\dot{\Phi}_j(s))\circ\Phi_j(s)ds)(x),$$  is Cauchy in the $L^\infty$ Hofer norm provided 
the sequence $ \Phi_j $ is Cauchy in the $C^0$ metric. Thus, 
the latter  converges in 
the complete metric space $\mathcal{N}^0([0,1]\times M\ ,\mathbb{R})$ to a time-independent continuous function denoted by $F_0$. On the 
other hand, the isotopy $$\Phi : t \mapsto \phi \circ \beta_t\circ \phi^{-1},$$ is a 
ssympeotopy which is a 1-parameter subgroup. Furthermore, it follows from the above statements that  the ssympeotopy
 $\Phi : t \mapsto \phi \circ \beta_t\circ \phi^{-1}$ is generated by $(F_0,\mathcal{K})$. Finally, we have constructed separately an 
example of 
 ssympeotopy which is a 1-parameter subgroup and whose generator is time independent. The strong flux of the 
ssympeotopy $\Phi : t \mapsto \phi \circ \beta_t\circ \phi^{-1}$ is given by the de Rham cohomology class of the harmonic 
 1-form $\mathcal{K}$, i.e. $\widetilde{Cal_0}(\Phi) = [\mathcal{K}].$  
The mass flow ssympeotopy $\Phi : t \mapsto \phi \circ \beta_t\circ \phi^{-1}$ is given by, 
$$ \widetilde{\mathfrak{F}}(\Phi)(h) =  
\dfrac{1}{(n-1)!}\int_M \mathcal{K}\wedge \omega^{n-1}\wedge h^\ast\sigma,
$$
for any  continuous function $h : M\rightarrow \mathbb S^1$. 
The Hofer-like length of the 
ssympeotopy $\Phi : t \mapsto \phi \circ \beta_t\circ \phi^{-1}$ is given by 
$l_\infty(\Phi) = osc(F_0) + |\mathcal{K}| = l_0(\Phi).$ 

\subsubsection*{Example 6.3}
Consider the torus $\mathbb T^{2n}$ with coordinates $(\theta_1,\dots,\theta_{2n})$ and 
equip it with the flat Riemannian metric $g_0 = \Sigma_{i = 1}^{2n} d\theta_i^2$. Then all the 1-forms $d\theta_i$, $i= 1,\dots, 2n$ are harmonics. 
Take the 1-forms $d\theta_i$ for $i= 1,\dots, 2n$ as basis for the space of harmonic 1-forms and consider the symplectic form
$
\omega = \sum_{i=1}^n d\theta_i \wedge d\theta_{i +n}.
$
\begin{itemize}
 \item Step (1). Given $v = (a_1,\dots, a_n, b_1,\dots, b_n) \in \mathbb R^{2n}$, 
the translation $x \mapsto x + v,$ on $ \mathbb R^{2n}$ induces a rotation $R_v$ on $\mathbb T ^{2n}$, 
which is a symplectic diffeomorphism.
Therefore, the smooth mapping $\{R_{v}^t\} : t \mapsto R_{tv}$ defines a harmonic isotopy 
generated by $(0, \mathcal{H})$ with $\mathcal{H} = \sum_{i=1}^n \left(a_i d\theta_{i+n} - b_i d\theta_i\right) .$ 
For each integer $j\textgreater 0$, consider the function 
$$f_j :[0,1]\rightarrow  \mathbb [0,1], t\mapsto \dfrac{j}{1+j}t$$ 
and denote by $\{R_v^{f_j}\}$ the smooth path $t \mapsto R_{v_j(t)}$ where $ v_j(t) = f_j(t)v.$
Since the sequence $(f_j)$ converges uniformly 
to the identity  mapping, we derive that $\{R_v^{f_j}\}\xrightarrow{\bar d} \{R_{v}^t\}$.   
For each $j$,  we have
$$\iota(\dot{R}^t_{v,j})\omega = \dfrac{j}{1+j}\sum_{i=1}^n \left(a_i d\theta_{i+n} - b_i d\theta_i\right).$$
Thus,
$$
l_\infty(\{R^t_{v,j}\}) = \dfrac{j}{1+j}|(-b_1,\dots, -b_n, a_1,\dots, a_n)| = l_0(\{R^t_{v,j}\}),
$$
where $ |(-b_1,\dots, -b_n, a_1,\dots, a_n)| = \sum_{i=1}^n |a_i|+ \sum_{i=1}^n |b_i| = |v|.$
 \item Step (2).  According to M\"{u}ller \cite{Mull08} the closed symplectic 
manifold $(\mathbb T^{2n}, \omega)$ admits 
at least a Hamiltonian homeomorphism which is not differentiable. Let us denote by $\phi$ such Hamiltonian homeomorphism. 
Then, by definition of $\phi$, 
there exists a sequence of Hamiltonian isotopies 
$ \Phi_j = \phi_{(U^j,0)}$ which is Cauchy in $(C^0 + L^\infty)$ such that $\phi_j := \Phi_j(1) \rightarrow \phi$ 
uniformly. Set, $$\eta_j^t = \dfrac{j}{1+j}\sum_{i=1}^n \left(a_i d\theta_{i+n} - b_i d\theta_i\right)
,$$ 
and 
$$ \mu^t(\eta_j,\Phi_j)  = \int_0^1\eta_j^t(\dot{\Phi}_j(s))\circ\Phi_j(s)ds,$$ 
for all $j$, and for all $t$. For each $j,$ the  symplectic isotopy 
$$\Psi_j : t \mapsto \phi_j^{-1} \circ R^t_{v,j}\circ \phi_j,$$ connects the identity to 
$ \phi_j^{-1} \circ R^1_{v,j}\circ \phi_j$. The latter  isotopy is generated by 
$( \mu(\eta_j,\Phi_j), \eta_j)$.

\item Step (3). By calculation we get, 
$$osc(\mu^t(\eta_j,\Phi_j) - \mu^t(\eta_{j +1},\Phi_{j +1}))\leq osc(\mu^t(\eta_j,\Phi_j) - \mu^t(\eta_{j},\Phi_{j +1})) $$
$$+ osc( \mu^t(\eta_{j},\Phi_{j +1}) - \mu^t(\eta_{j +1},\Phi_{j +1})),$$
for all $t$, and for all $j$. Thus, Lemma 3.9 found in \cite{TD} implies that,
$$ \max_t osc(\mu^t(\eta_j,\Phi_j) - \mu^t(\eta_{j},\Phi_{j +1}))\rightarrow0, j\rightarrow\infty,$$ 
while Lemma 3.4 found in \cite{BanTchu} implies that 
$$\max_t osc( \mu^t(\eta_{j},\Phi_{j +1}) - \mu^t(\eta_{j +1},\Phi_{j +1}))\rightarrow0, j\rightarrow\infty.$$ 
So, the sequence $\mu(\eta_j,\Phi_j)$ 
is Cauchy in the $L^\infty$ Hofer norm, and then converges in 
the complete metric space $\mathcal{N}^0([0,1]\times\mathbb T^{2n} ,\mathbb{R})$ to a continuous family 
of continuous functions denoted by $F_1$. Thus, 
$$( \Psi_j, (\mu(\eta_j,\Phi_j) , \eta_j) ) \xrightarrow{L^\infty + C^0}
( \phi^{-1} \circ R_{v}\circ \phi, (F_1 , \sum_{i=1}^n \left(a_i d\theta_{i+n} - b_i d\theta_i\right)).$$ Note that the path  
 $(t,x) \mapsto (\phi^{-1} \circ R_{v}^t\circ \phi)(x),$ is continuous, and not differentiable.

\item Step (4). The paths $\Psi : (t,x) \mapsto (\phi^{-1} \circ R_{v}^t\circ \phi)(x)$ just constructed above 
yields a non-trivial ssympeotopy which is not differentiable. The strong flux of the 
ssympeotopy $\Psi$ is given by the de Rham cohomology class of the harmonic 
 1-form $ \sum_{i=1}^n \left(a_i d\theta_{i+n} - b_i d\theta_i\right)$, i.e. 
$$\widetilde{Cal_0}(\Psi) = [\sum_{i=1}^n \left(a_i d\theta_{i+n} - b_i d\theta_i\right)].$$  
The mass flow ssympeotopy $\Psi$ is given by, 

$$ \widetilde{\mathfrak{F}}(\Psi)(h) = 
\dfrac{1}{(n-1)!}\sum_{i=1}^n\int_{\mathbb T^{2n}}\left(a_i d\theta_{i+n} - b_i d\theta_i\right)\wedge \omega^{n-1}\wedge h^\ast\sigma,
$$
for any continuous function $h : \mathbb T^{2n}\rightarrow \mathbb S^1$.\\ The Hofer-like length of the 
ssympeotopy $\Psi$ is given by 
$$l_\infty(\Psi) =  \max_t osc(F_1^t) + |v|.$$
 
\end{itemize}

\subsubsection*{Example 6.4} Similarly as in Example 6.3, consider the torus $\mathbb T^2$ as the square:
$$
 \{(p,q) \mid 0\leq p \leq 1 \text{, } 0\leq q \leq 1\} \subset \mathbb R^2
$$
with opposite sides identified. Let $\mathbb D^2\subset \mathbb R^2$ be the $2-$disk of radius $\tau \in ]0, 1/8[$ centered at 
$A = (a, 0)$ with 
$7/8\leq a \textless 1$, and  let $ \Lambda^2(\tau)$ be the corresponding subset in $ \mathbb T^2$. 
 For any $\nu \textless 1/4$, consider 
$$
B(\nu) = \{(x,y) \mid 0\leq x \textless \nu\} \subset \mathbb{R}^2,
$$
and let $ C(\nu)$ be the corresponding subset in $ \mathbb T^2$. 
Let $\phi$ be a Hamiltonian homeomorphism of $ \mathbb T^2$ supported in  $\Lambda^2(\tau)$. 
If $v = (a_1,0)$ with $\nu\leq a_1 \leq 1/2$,
then the map $\phi^{-1} \circ R_{v}\circ \phi$ displaces completely the set $C(\nu)$, where $R_{v}$ 
is the rotation induced by the translation $(p,q) \mapsto (p + a_1, q)$. 
It follows from Example 6.3 that, 
$$
0\textless e_S(C(\nu))\leq e^0_S( C(\nu) )\leq \|\phi^{-1} \circ R_{v}\circ \phi\|_{SHL}$$
$$ \leq \dfrac{\max_t osc(F_0^t) + \max_t osc(G_0^t)}{2} +  \dfrac{1}{2} \textless \infty, 
$$
where $F_0$ and $G_0$ are two continuous families of continuous functions. Thus, 
 the strong symplectic displacement energy $e^0_S( C(\nu) ) $ of the set $ C(\nu)$ is finite and positive.

\begin{center}
 {\bf Acknowledgments:}
 \end{center}
\begin{center}
I thank the Professors Barry Green and Jeff Sanders for awarding to me a visiting researcher fellowships at 
AIMS South Africa. This stay was helpful in the finalization of this note.\\
\end{center}
\begin{center}
I thank Peter Spaeth and Stefan M\"{u}ller for helpful hints.\\ 
\end{center}

\begin{center}
I would also thank the following Professors for encouragements : Idris Assani, Diallo Assimiou, Cyriaque Atindogbe, David B\'ekoll\`e, N\"{o}el Lohoue, 
Ferdinand Ngakeu, Boyom Nguiffo, Boniface Nkemzi, 
  Carlos Ogouyandjou, Eugene OKassa, Anatole Temgoua, Leonard Todihounde, J\"{o}el Tossa, Aissa Wade, Fran\c{c}ois Wamon.\\
\end{center}

\bibliographystyle{plain}

\begin{center}

\end{center}

\begin{center}
 $\ast$
Department of Mathematics\\ 
The University of Buea,\\ South West Region, Cameroon
\end{center}

\end{document}